# A RATE OF CONVERGENCE RESULT FOR THE LARGEST EIGENVALUE OF COMPLEX WHITE WISHART MATRICES[1]

BY NOUREDDINE EL KAROUI

*University of California, Berkeley*

It has been recently shown that if $X$ is an $n \times N$ matrix whose entries are i.i.d. standard complex Gaussian and $l_1$ is the largest eigenvalue of $X^*X$, there exist sequences $m_{n,N}$ and $s_{n,N}$ such that $(l_1 - m_{n,N})/s_{n,N}$ converges in distribution to $W_2$, the Tracy–Widom law appearing in the study of the Gaussian unitary ensemble. This probability law has a density which is known and computable. The cumulative distribution function of $W_2$ is denoted $F_2$.

In this paper we show that, under the assumption that $n/N \to \gamma \in (0, \infty)$, we can find a function $M$, continuous and nonincreasing, and sequences $\tilde{\mu}_{n,N}$ and $\tilde{\sigma}_{n,N}$ such that, for all real $s_0$, there exists an integer $N(s_0, \gamma)$ for which, if $(n \wedge N) \geq N(s_0, \gamma)$, we have, with $l_{n,N} = (l_1 - \tilde{\mu}_{n,N})/\tilde{\sigma}_{n,N}$,

$$\forall s \geq s_0 \qquad (n \wedge N)^{2/3} |P(l_{n,N} \leq s) - F_2(s)| \leq M(s_0) \exp(-s).$$

The surprisingly good 2/3 rate and qualitative properties of the bounding function help explain the fact that the limiting distribution $W_2$ is a good approximation to the empirical distribution of $l_{n,N}$ in simulations, an important fact from the point of view of (e.g., statistical) applications.

**1. Introduction.** Very important progress has been made in recent years in our understanding of the behavior of the eigenvalues of a large number of large-dimensional random matrices. Many new results concern the fluctuation of these eigenvalues: we now have convergence in distribution results, as opposed to maybe more classical almost-sure convergence statements. These new findings show great promise for applications, in particular, in Statistics.

Received October 2004; revised December 2005.
[1]Supported in part by Grants NSF DMS-00-77621 and ANI-00-8584 (ITR).

Our focus in this paper is on fine convergence properties of the largest eigenvalue of a class of random covariance matrices. Our work owes a lot to the pioneering work of Tracy and Widom [20, 21, 22] and also to that of Johnstone [11]. Exciting recent developments in the area of random covariance matrices and/or Tracy–Widom distributions can also be found, for instance, in [2, 3, 4, 5, 6].

Let us now be more specific about the question we address. In a series of papers by Forrester [8], Johansson [10] and Johnstone [11], it was shown that if $X$ is an $n \times N$ matrix with i.i.d. standard complex Gaussian entries, $l_1$, the largest eigenvalue of $X^*X$, when properly renormalized, converges in distribution [when $n/N \to \gamma \in (0,\infty)$ as $n \to \infty$] to the Tracy–Widom law appearing in the study of the Gaussian unitary ensemble, which we call $W_2$. The paper [7] shows that this result holds even when $n/N \to \infty$ or 0. We recall that a random variable $Z$ is said to have a standard complex Gaussian distribution [denoted $\mathcal{N}_{\mathbb{C}}(0,1)$] if $Z = (Z_1 + iZ_2)/\sqrt{2}$, where $Z_1$ and $Z_2$ are i.i.d. $\mathcal{N}(0,1)$.

The chronology of the problem is the following: Forrester worked on the case $n - N$ constant (see Section 3.5 in [8]), Johansson (Theorem 1.6 in [10]) showed that the result held when $n = \gamma N + a_N$, with $a_N = \mathrm{O}(N^{1/3})$, and Johnstone (Section 3 and Appendix A.6 in [11]), while being primarily interested in the case of real entries, relaxed the assumption to $n/N \to \gamma \in (0,\infty)$. In [7] we show that the result holds even when $n/N \to 0$ or $\infty$.

We denote by $F_2$ the cumulative distribution function of $W_2$. Following [11], we call matrices like $X^*X$, that is, complex Wishart matrices with identity covariance, complex white Wishart matrices. In random matrix theory terminology, we are therefore considering the Laguerre unitary ensemble (LUE). We remind the reader (see [20]) that if $q$ solves

$$q''(x) = xq(x) + 2q^3(x) \quad \text{and} \quad q(x) \sim \mathrm{Ai}(x) \qquad \text{as } x \to \infty,$$

then

$$F_2(s) = \exp\biggl(-\int_s^\infty (x-s)q^2(x)\,dx\biggr).$$

Let us recall the main result we are using from [8, 10] and [11]. Let

$$n_+ = n + 1/2 \quad \text{and} \quad N_+ = N + 1/2,$$

$$\mu_{n,N} = (\sqrt{n_+} + \sqrt{N_+})^2 \quad \text{and} \quad \sigma_{n,N} = (\sqrt{n_+} + \sqrt{N_+})(n_+^{-1/2} + N_+^{-1/2})^{1/3}.$$

The following was progressively shown in these papers:

THEOREM 1. *If $X$ is an $n \times N$ matrix with i.i.d. $\mathcal{N}_{\mathbb{C}}(0,1)$ entries and $n/N \to \gamma \in (0,\infty)$ when $n \to \infty$,*

$$\frac{l_1(X^*X) - \mu_{n,N}}{\sigma_{n,N}} \Longrightarrow W_2.$$

To be completely correct, we need to say that the results in [8, 10, 11] are not stated exactly as above, but they imply their respective version of Theorem 1 in the aforementioned form. In particular, the different authors of [8, 10, 11] use different centering and scaling sequences. For instance, in [11], the real analog of Theorem 1 is stated with $n_+$ replaced by $n-1$ and $N_+$ replaced by $N$. However, $\mu_{n,N}$ and $\sigma_{n,N}$ are the natural sequences from the point of view of the analysis presented in [11] [see Section 3 and equations (A.7) and (A.8) there].

It is clear that to use the Tracy–Widom law in Statistics or other applied areas, we would like to know what is the rate of convergence in Theorem 1. In other words, we would like to have a result similar to the Berry–Esseen refinement of the central limit theorem. The celebrated Berry–Esseen theorem states that the error made by approximating the distribution of the mean of a sample of size $m$ of i.i.d. random variables with a third moment by a Gaussian distribution is of size $m^{-1/2}$. Here we will show that, for the properly renormalized largest eigenvalue of $X^*X$, a similar measure of error is of size at most $m^{-2/3}$, where now $m = \min(n, N)$. We will also shed light on the involved issue of choosing good recentering and rescaling sequences, which significantly improve the speed of convergence result.

We need to introduce some notation before stating the theorem. We call our final centering and scaling sequences $\tilde{\mu}_{n,N}$ and $\tilde{\sigma}_{n,N}$, and define them below. We will denote by

$$l_{n,N} = \frac{l_1(X^*X) - \tilde{\mu}_{n,N}}{\tilde{\sigma}_{n,N}}$$

the sequence of centered and scaled largest eigenvalues.

THEOREM 2. *Let us assume that $N \leq n$, and that $n/N \to \gamma \in [1, \infty)$. We call*

$$\gamma_{n,N} = \frac{\mu_{n-1,N}\sigma_{n,N-1}^{1/2}}{\mu_{n,N-1}\sigma_{n-1,N}^{1/2}},$$

$$\tilde{\sigma}_{n,N} = (1 + \gamma_{n,N})\left(\frac{1}{\sigma_{n-1,N}} + \frac{\gamma_{n,N}}{\sigma_{n,N-1}}\right)^{-1}$$

*and*

$$\tilde{\mu}_{n,N} = \left(\frac{1}{\sigma_{n-1,N}^{1/2}} + \frac{1}{\sigma_{n,N-1}^{1/2}}\right)\left(\frac{1}{\mu_{n-1,N}\sigma_{n-1,N}^{1/2}} + \frac{1}{\mu_{n,N-1}\sigma_{n,N-1}^{1/2}}\right)^{-1}.$$

*Let $l_{n,N}$ denote the largest eigenvalue of $X^*X$ renormalized as above.*

*There exists a function $M$, such that, for all real $s_0$, there exists an integer $N(s_0, \gamma)$ for which we have*

$$\forall s \geq s_0, \forall N \geq N(s_0, \gamma) \qquad N^{2/3}|P(l_{n,N} \leq s) - F_2(s)| \leq M(s_0)e^{-s}.$$

*If $N > n$, the theorem is valid after we switch $n$ and $N$ in all the displays. This effectively shows that it holds for $\gamma \in \mathbb{R}_+^*$. The function $M$ can be chosen to be continuous and nonincreasing. It has a finite limit at $+\infty$, and is unbounded at $-\infty$.*

*The theorem remains valid under a wider set of technical conditions on $\tilde{\mu}_{n,N}$ and $\tilde{\sigma}_{n,N}$ detailed in Section* 4.1. *More information about $M$ can be found in Section* 4.2.

To explain the importance of the centering and scaling sequences, we have to say that with $\mu_{n,N}$ and $\sigma_{n,N}$ as renormalization sequences, we could not show that the speed of convergence was higher than $N^{1/3}$. The proof of Theorem 2 makes clear that $\mu_{n,N}$—the centering sequence—was probably in this case the quantity that was "slowing us down."

It should be noted that the bounding function $M$ and other bounding functions we will obtain in the course of the proof look like they might depend on the particular $\gamma$ we are choosing. We will see in Appendix A.6 that we can actually bound them independently of $\gamma$.

Before we proceed, we wish to mention that Iain Johnstone, while working on an analog of Theorem 1 for the so-called Jacobi ensemble of random matrices [12], found that the pointwise convergence of the corresponding kernel to the Airy kernel could happen at rate 2/3 (see Section 2 for more information on how eigenvalue problems are turned into integral operator ones, with explicit kernels). He then adapted his computations to the Laguerre case we are considering and found the same pointwise result. Though such pointwise techniques do not yield similar estimates for convergence in trace class norm, which is the delicate notion of convergence we have to use for the problems under investigation, his computations suggested that the 2/3 rate in Theorem 2 was a plausible goal. We wish to point out that an important piece of our analysis, Lemma 2, highlights a phenomenon that is similar in spirit—but different in substance—to the one observed by Iain Johnstone in his pointwise computations. We are very grateful to him for sharing his insights with us.

The plan of the article is as follows: in Section 2 we review the techniques used in [11] to obtain Theorem 1 and outline our strategy for the proof. We then proceed to give the proof (in Section 3) of Theorem 2 in two steps. We will first show a naive analysis, with nonrefined centering and scaling sequences, to get intermediate (and needed) results and to highlight the difficulty that arises. We then provide a solution that makes natural our choice of $\tilde{\mu}_{n,N}$ and $\tilde{\sigma}_{n,N}$. Finally, we discuss in Section 4 some of the properties of the bounding function $M$ and present simulations to assess the quality of the Tracy–Widom approximation—using our centering and scaling sequences—across a range of dimensions.

Most of the technical questions that are necessary to carry out the proof but would obscure its explanation are relegated to the appendices. In the rest of the paper we will assume that $n \geq N$, $n/N \to \gamma \in [1, \infty)$. Since $X$ and $X^*$ have the same largest singular value, this will take care of $n < N$ and $\gamma \in (0, 1]$.

## 2. Outline of the proof.

2.1. *Review of known results.* A strength of the method used in [11] is that the intermediate steps lead to finite-dimensional equalities, and limits are taken only at the last step. This is a crucial element in our being able to get rates of convergence estimates.

Recall that Johnstone, using Tracy and Widom's work [21, 22], shows that

$$\mathbf{P}\left(\frac{l_1 - \mu_{n,N}}{\sigma_{n,N}} \leq s\right) = \det(\mathrm{Id} - S_\tau),$$

where $S_\tau$ is a known kernel acting here on $L^2([s, \infty))$, and det is understood as a Fredholm determinant. As an aside, let us remark that the previous equation holds for arbitrary choices of centering and scaling sequences. Naturally, the resulting operator depends strongly on which choice is made, but for all choices such an operator exists. In particular, in [11], the operator $S_\tau$ (and, hence, the index $\tau$) depends on the sequences $\mu_{n,N}$ and $\sigma_{n,N}$ (see Section 3 and equations (A.7) and (A.8) in [11]). We will nevertheless use the notation $S_\tau$ even when working with different sequences $\tilde{\mu}_{n,N}$ and $\tilde{\sigma}_{n,N}$ (instead of $\mu_{n,N}$ and $\sigma_{n,N}$), in order to avoid cumbersome notation like $S_{\tilde{\tau}}$, but we will keep in mind that our $S_\tau$ is not necessarily exactly the one appearing in [11]. We now go back to our argument. Using continuity of the Fredholm determinant with respect to trace class norm, Johnstone proceeds to show that $S_\tau \to \bar{S}$ in trace class norm, where $\bar{S}$, the Airy kernel, was shown in [20] to satisfy

$$F_2(s) = \det(\mathrm{Id} - \bar{S}).$$

It turns out that we can even control the difference of two Fredholm determinants in this situation using the following result from [18], cited in [17] (Lemma 4, pages 323–324) and [9] (Section II.4 and Theorem IV.5.2):

LEMMA 1 (Seiler–Simon). *Let $A$ and $B$ be in $\mathcal{S}_1$, the family of all trace class operators on a separable Hilbert space $\mathcal{H}$. If $\|\cdot\|_1$ represents trace class norm, we have*

$$|\det(\mathrm{Id} - A) - \det(\mathrm{Id} - B)| \leq \|A - B\|_1 \exp(\|A\|_1 + \|B\|_1 + 1).$$

2.2. *Strategy for the proof.* The strategy is now clear. Since all we want to control is the difference $|\det(\mathrm{Id} - S_\tau) - \det(\mathrm{Id} - \bar{S})|$, the previous display makes clear that we just need to study $\|S_\tau - \bar{S}\|_1$.

As a matter of fact, $\|\bar{S}\|_1(s)$ is, when seeing $\bar{S}$ as an operator on $L^2([s,\infty))$, a well understood function of $s$ and will not cause any problems. We give more detail in Section 4.2. Working toward the control of $\|S_\tau - \bar{S}\|_1$, we recall that we know from [11] that

$$S_\tau = H_\tau G_\tau + G_\tau H_\tau$$

and

$$\bar{S} = 2G^2,$$

for explicit $G$, $G_\tau$ and $H_\tau$, all of them being Hilbert–Schmidt operators on $L^2([s,\infty))$, for all $s \in \mathbb{R}$.

We have the following elementary lemma:

LEMMA 2.

$$2\|S_\tau - \bar{S}\|_1 \leq \|G_\tau + H_\tau - 2G\|_2 \|G_\tau + H_\tau + 2G\|_2 + \|G_\tau - H_\tau\|_2^2,$$

*where $\|\cdot\|_2$ represents the Hilbert–Schmidt norm on $L^2([s,\infty))$.*

Since the tools involved in the proof of this lemma are somewhat different from the ones involved in other proofs, this is proved in Appendix A.4.

So to prove Theorem 2, it will be sufficient to show the following:

LEMMA 3. *Assume that $n/N \to \gamma \in [1,\infty)$. Then we have $\forall s_0 \in \mathbb{R}$, $\exists N(s_0, \gamma), \forall s \geq s_0, \forall N > N(s_0, \gamma)$,*

(P1) $\qquad N^{2/3}\|G_\tau + H_\tau - 2G\|_2 \leq C(s_0)\exp(-s/2),$

(P2) $\qquad N^{1/3}\|G_\tau - G\|_2 \leq C(s_0)\exp(-s/2),$

(P3) $\qquad N^{1/3}\|H_\tau - G\|_2 \leq C(s_0)\exp(-s/2),$

*where $C$ is a continuous, nonincreasing function.*

The rest of the article will be devoted to proving these estimates. Let us now make a few remarks. The first one is structural: our proof makes fundamental use of the structure of $S_\tau$, that is, the fact that it is the sum of the product of two Hilbert–Schmidt operators. It also heavily relies on the fact that those operators [seen as acting on $L^2([s,\infty))$] are kernel operators with kernels of the form $K_s(x,y) = K(x+y-s)$, which reduces our problem to studying certain functions, as opposed to the potentially more complicated objects that are general operators. Other problems having this same structure could be attacked by the same approach.

Closer to the problem we actually work on, let us mention that we have the choice of which rate to increase from $1/3$ to $2/3$ in Lemma 3. In our decomposition of $S_\tau - \bar{S}$, $G_\tau + H_\tau - 2G$ plays a more important role (in terms of rates) than $G_\tau - G$ or $H_\tau - G$ do. This is why an important effort will be devoted to showing that we can get rate $2/3$ for the convergence of $G_\tau + H_\tau - 2G$ in the Hilbert–Schmidt norm to 0. We will show later that the "natural" centering and scaling (natural for the perturbation analysis) leads to rate $2/3$ for $\|G_\tau - G\|_2$, but rate $1/3$ for the two other elements of the previous display. Our "optimized" centering reflects the fact that we had to find a trade-off between an optimal centering for $G_\tau$ and an optimal centering for $H_\tau$. This is partly why the "optimized" sequences look so involved.

In other respects, the operators mentioned here are kernel operators whose kernels are well known and understood. Because these kernels can be related to the solution of a perturbed Airy equation, our task is now essentially reduced to studying in detail the properties of a solution of a certain differential equation. This will become more clear in the course of the article. Note, finally, that Lemma 3 only deals with Hilbert–Schmidt norms, which are considerably simpler to manipulate and bound than trace class norms.

2.2.1. *Technical details on the elements of the problem.* Because of elementary results of linear algebra, namely, the fact that $X$ and $X^*$ have the same singular values (except for the multiplicity of those at 0), we can and do assume that $n \geq N$ in what follows. The notation we use are the ones found in [11].

Let us introduce $\alpha_N = n - N$, and

$$\varphi_k(x; \alpha) = \sqrt{\frac{k!}{(k+\alpha)!}} x^{\alpha/2} e^{-x/2} L_k^\alpha(x),$$

where $L_k^\alpha$ is the $k$th Laguerre polynomial associated with $\alpha$. Then calling $a_N = \sqrt{Nn}$, $\varphi$ and $\psi$ are defined as

$$\varphi(x) = (-1)^N \sqrt{\frac{a_N}{2}} \varphi_N(x; \alpha_N - 1) x^{-1/2},$$

$$\psi(x) = (-1)^{N-1} \sqrt{\frac{a_N}{2}} \varphi_{N-1}(x; \alpha_N + 1) x^{-1/2}.$$

Now let Ai denote the Airy function, $\varphi_\tau(s) = \tilde{\sigma}_{n,N} \varphi(\tilde{\mu}_{n,N} + s\tilde{\sigma}_{n,N})$ and similarly $\psi_\tau(s) = \tilde{\sigma}_{n,N} \psi(\tilde{\mu}_{n,N} + s\tilde{\sigma}_{n,N})$. Finally, we have, when considering $G_\tau$, $H_\tau$ and $G$ as operators on $L^2([s, \infty))$,

$$G_\tau(x,y) = \varphi_\tau(x+y-s), \qquad G(x,y) = 2^{-1/2} \text{Ai}(x+y-s),$$
$$H_\tau(x,y) = \psi_\tau(x+y-s).$$

The $\varphi$ and $\psi$ we introduce here are the same as in [11], but we have used some elementary properties of Laguerre polynomials (see [19], page 102) to simplify their expression. All the details are given in Appendix A.5, where one will also find a remark explaining why the case $\alpha_N = 0$ does not pose a problem. This makes our rate work much simpler later, since we will have to deal with only two pieces (instead of four if we had kept the original representation) when adjusting the centering and scaling sequences to get rate $2/3$.

2.2.2. *Formulating the problem.* In his study, Johnstone [11] introduced an "intermediate" function that proved to be a powerful tool in the study of the problem. We introduce, following his idea, the very closely related function

$$F_{n,N}(z) = (-1)^N \sigma_{n,N}^{-1/2} \sqrt{N!/n!} z^{(\alpha_N+1)/2} e^{-z/2} L_N^{\alpha_N}(z). \tag{1}$$

This is valuable as the study of $F_{n,N}$ as $N$ gets large turns out to be the study of a perturbation of the Airy equation, and this has been the subject of extensive investigations. See, for instance, [15], Chapter 11. In special functions parlance, we have

$$F_{n,N}(z) = \sigma_{n,N}^{-1/2} \frac{1}{\sqrt{n!N!}} W_{\kappa_N, \lambda_N}(z),$$

where $W_{\kappa_N, \lambda_N}$ stands for the Whittaker function with parameters $\kappa_N = N + (\alpha_N + 1)/2$ and $\lambda_N = \alpha_N/2$.

In the situation we are investigating here, namely, $n/N$ tends to a finite limit as $n$ and $N$ get large, [11] introduces the investigations of Olver [15] concerning the so-called differential equations with one turning point to study the properties of $W_{\kappa_N, \lambda_N}$. In the case $n/N \to \infty$, we would have to take care of the fact that the turning points coalesce (see [7]) to do the same. Note that we cannot really apply Olver's results in [15] directly, since the differential equation we consider depends on two parameters. In Appendix A.3, we explain why they apply nonetheless. This was implicit in [11], but we make it explicit for the sake of completeness.

We now need to relate $F_{n,N}$ to $\varphi(\cdot; \alpha)$ and $\varphi_\tau$. We remark that

$$\varphi_N(x; \alpha_N) = (-1)^N \sigma_{n,N}^{1/2} x^{-1/2} F_{n,N}(x),$$

$$\varphi(x) = \sqrt{\frac{a_N}{2}} \sigma_{n-1,N}^{1/2} F_{n-1,N}(x)/x$$

and

$$\varphi_\tau(x) = \frac{1}{\sqrt{2}} \left( \frac{\sqrt{a_N} \sigma_{n-1,N}^{1/2} \tilde{\sigma}_{n,N}}{\mu_{n-1,N}} \right) F_{n-1,N}(\tilde{\mu}_{n,N} + \tilde{\sigma}_{n,N} x) \left( \frac{\mu_{n-1,N}}{\tilde{\mu}_{n,N} + \tilde{\sigma}_{n,N} x} \right).$$

We have similar expressions for $\psi_\tau$, if we replace $(n-1, N)$ by $(n, N-1)$ in the previous expression.

So the problem essentially reduces to having a good understanding of $F_{n-1,N}$ and the other quantities written in the equation defining $\varphi_\tau$.

2.2.3. *What do we need to control in order to get the rates?* Our objective is to get the rates mentioned in Lemma 3. In order to do this, we will prove the following fact:

FACT 2.2.1. *Suppose $n/N \to \gamma \in [1, \infty)$. Let us call $G(x) = \text{Ai}(x)/\sqrt{2}$. We have $\forall s_0 \exists N(s_0, \gamma) \forall z \in [s_0, \infty), \forall N > N(s_0, \gamma)$,*

$$N^{2/3} |\varphi_\tau(z) + \psi_\tau(z) - 2G(z)| \leq C(s_0) \exp(-z/2), \tag{2}$$

$$N^{1/3} |\varphi_\tau(z) - G(z)| \leq C(s_0) \exp(-z/2), \tag{3}$$

$$N^{1/3} |\psi_\tau(z) - G(z)| \leq C(s_0) \exp(-z/2), \tag{4}$$

*for a function $C$ that is continuous and nonincreasing on $\mathbb{R}$.*

The reason this implies Lemma 3 is that all the kernels—when we see the operators as acting on $L^2([s, \infty))$—are of the form $K_s(x, y) = K(x + y - s)$, and we are dealing with Hilbert–Schmidt norms. In somewhat more detail, let us consider a Hilbert–Schmidt operator $A$ [on $L^2([s, \infty))$] with kernel $K_s(x, y)$. Let us assume that $K_s(x, y) = K(x + y - s)$ and $n^\beta |K(z)| \leq \chi \exp(-z/2)$, for a certain $\chi \in \mathbb{R}_+$. Let us call $\|\cdot\|_2(s)$ the Hilbert–Schmidt norm on $L^2([s, \infty))$. It is well known (see [16], Theorem VI.23, page 210) that $\|A\|_2^2(s) = \iint_{[s,\infty)^2} |K_s(x,y)|^2 \, dx \, dy$. So in our situation, since $|K_s(x,y)|^2 = |K(x+y-s)|^2 \leq n^{-2\beta} \chi^2 \exp(-(x+y-s))$, we immediately have $\|A\|_2(s) \leq n^{-\beta} \chi \exp(-s/2)$. So Fact 2.2.1 clearly implies Lemma 3.

**3. Proof.** Before we proceed to giving the proof, we remind the reader of the analysis carried out in [11]: $F_{n,N}$ (actually a closely related quantity) was analyzed using the Liouville–Green approximation, which we detail now. Pushing further this type of analysis will give us explicit bounds on how far $F_{n,N}$ deviates from the Airy function and will be the centerpiece of the analysis. A remark about notation: they are "naturally" heavy and we warn the reader that some abuse will take place, as we will not always use both indices "$n$" and "$N$," even when it is clear that the function or sequence depends on both of them. But they will always be present when we give rigorous arguments and there might be a doubt about the quantities we are talking about.

This section is organized as follows: we first recall in detail the method that was used in [11] to show convergence of $S_\tau$ to $\bar{S}$. We then prove an intermediary result concerning the rate of convergence $(2/3)$ of $F_{n,N}(\mu_{n,N} + \sigma_{n,N}x)(\mu_{n,N}/(\mu_{n,N}+\sigma_{n,N}x))$ to $\operatorname{Ai}(x)$ on semi-infinite intervals of the form $[s_0, \infty)$. Then, we will show (in Section 3.4) that we can only guarantee a rate of $1/3$ if we slightly perturb the centering and scaling. This will show equations (P2) and (P3) of Lemma 3. Finally, we will prove equation (P1) of Lemma 3 in Section 3.5.

3.1. *Liouville–Green approximation for Whittaker functions.* We recall that if we call $\alpha_N = n - N$, $w_N(x) = x^{(\alpha_N+1)/2}e^{-x/2}L_N^{\alpha_N}(x)$, we have

(Whittaker) $$\frac{d^2w_N}{dx^2} = \left\{\frac{1}{4} - \frac{\kappa_N}{x} + \frac{\lambda_N^2 - 1/4}{x^2}\right\}w_N,$$

with $\kappa_N = N + \frac{\alpha_N+1}{2}$ and $\lambda_N = \frac{\alpha_N}{2}$.

Following [11], after changing variables to $\xi = x/\kappa_N$, we get

$$\frac{d^2w_N}{d\xi^2} = \{\kappa_N^2 f(\xi) + g(\xi)\}w_N,$$

with

$$f(\xi) = \frac{(\xi - \xi_1)(\xi - \xi_2)}{4\xi^2} \quad \text{and} \quad g(\xi) = -\frac{1}{4\xi^2}.$$

Here $\xi_1 = 2 - \sqrt{4 - \omega_N^2}$ and $\xi_2 = 2 + \sqrt{4 - \omega_N^2}$ with $\omega_N = \frac{2\lambda_N}{\kappa_N} = 2\frac{n-N}{n+N+1}$. Under our assumptions about $n$ and $N$, $\omega_N^2 \in [0, 4 - \delta]$, with $\delta > 0$. This prevents the turning points $\xi_1$ and $\xi_2$ from coalescing, a crucial point in what follows. Note that, to be precise, we should write that $f(\xi)$ is really $f(\omega_N, \xi)$, as the family of functions $f$ is indexed by $\omega_N$.

Now, following [15], Chapter 11, [11] introduces the change of variable (known as the Liouville–Green or WKB method)

$$\tfrac{2}{3}\zeta^{3/2} = \int_{\xi_2}^{\xi} f^{1/2}(t)\, dt.$$

Defining a new dependent variable $W$ by $w = (d\zeta/d\xi)^{-1/2}W$, one gets the new differential equation

$$\frac{d^2W}{d\zeta^2} = \{\kappa_N^2\zeta + \upsilon(\omega_N, \zeta)\}W,$$

where we insist on the fact that $\upsilon$ is also a function of $\omega$, because $f$ was. We also note that $\upsilon$ has a somewhat explicit expression, to which we come back in Appendix A.3, so the reader is referred there for more details. Then, if we denote by $\hat{f} = f/\zeta$, the recessive solution of the equation (Whittaker) satisfies (see [11], equation (5.16))

$$w_N(\kappa_N\xi) = c_{\kappa_N}\hat{f}^{-1/4}(\xi)\{\operatorname{Ai}(\kappa_N^{2/3}\zeta) + \varepsilon_2(\kappa_N, \xi)\},$$

where

$$|\varepsilon_2(\kappa_N,\xi)| \le \mathbf{M}(\kappa_N^{2/3}\zeta)\mathbf{E}^{-1}(\kappa_N^{2/3}\zeta)\left[\exp\left(\frac{\lambda_0}{\kappa_N}F(\omega_N)\right) - 1\right].$$

**M** and **E** are the modulus and weight functions introduced in [15], Chapter 11, pages 394–397. Also, $\lambda_0 \doteq 1.04$, and $F(\omega)$ is bounded when $\omega^2$ is in $[0, 4-\delta]$ (see Appendix A.3). The sequence $c_{\kappa_N}$ has an explicit expression, and a little algebra (see also [11], equation (5.20)) leads to

$$(5) \qquad F_{n,N}(x) = r_N\left(\frac{\kappa_N}{\sigma_{n,N}^3}\right)^{1/6}\hat{f}^{-1/4}(\xi)\{\text{Ai}(\kappa_N^{2/3}\zeta) + \varepsilon_2(\kappa_N,\xi)\},$$

with $x = \kappa_N\xi$. If we denote by $n_+ = n + 1/2$ and $N_+ = N + 1/2$,

$$r_N^2 = \frac{2\pi\exp[-(n_+ + N_+)]n_+^{n_+}N_+^{N_+}}{N!n!},$$

with $r_N$ being nonnegative (as an aside, note that our $r_N$ corresponds to [11]'s $(-1)^N r_N$). We will occasionally use only the index $N$ instead of $n, N$. When needed, we will make precise the $n$ we are dealing with.

3.2. *Gathering the different elements*: *a useful intermediate result.* This subsection sets up the core of the technical work of the article. We introduce a class of functions (containing, of course, $\varphi_\tau$ and $\psi_\tau$) and study in Section 3.3 how its members deviate from the Airy function; there, we are going to choose the equivalent of an optimal centering and scaling for the $(n, N)$ pair and show that we can obtain rate $2/3$ [recall that $\varphi$ corresponds to $(n-1, N)$ and $\psi$ corresponds to $(n, N-1)$]. At the end of Section 3.3 we will have a much finer understanding of the issues involved, will know what is easy and hard to deal with and will have all the technical elements needed to tackle the proof of equation (P1). The conclusion of our work is that one could actually have rate $2/3$ in (P2) and (3). The proofs will also imply that getting the rate $1/3$ is "easy," and so we should reserve the hard-earned $2/3$ rate to something that is harder to deal with, namely, equation (P1). Section 3.3 will give us all the elements needed to do this, the details being taken care of in Section 3.5. We now carry out in detail the analysis.

3.2.1. *About the deviation of $F_{n,N}$ and related functions from* Ai. Let us note that $\varphi_\tau(s)$ and $\psi_\tau(s)$ have the same functional form. They can be written

$$\frac{1}{\sqrt{2}}\alpha_{j,k}F_{j,k}(x(s))\left(\frac{\mu_{j,k}}{x(s)}\right),$$

where $j$ and $k$ are integer indexes, $F_{j,k}$ stands for the function introduced in equation (1), $\mu_{j,k}$ is the sequence mentioned in the Introduction, $x(s) =$ $m_{j,k} + s_{j,k}s$ and $\alpha_{j,k}$, $m_{j,k}$ and $s_{j,k}$ are sequences of numbers independent of $s$. Getting a lower bound on the rate of convergence in our original problem essentially reduces to studying how members of this class of functions deviate from the Airy function. A central element in doing this is the expression for $F_{n,N}$ mentioned in equation (5). We will decompose it into four blocks: $r_N$, $(\kappa_N/\sigma_{n,N}^3)^{1/6}\hat{f}^{-1/4}$, $\text{Ai}(\kappa_N^{2/3}\zeta)$ and $\varepsilon_2(\kappa_N,\xi)$.

The easiest part to deal with is $\varepsilon_2$, as we will see, it is controlled by $1/\kappa_N$, with $\kappa_N \sim (1+\gamma)N/2$, and our aim is to get only speed $N^{2/3}$. So it is not going to cause us any problems (we will show that we control the other part of the expression defining $\varepsilon_2$).

Given the expression immediately following equation (5), $r_N$ has easy to analyze asymptotics, and as is shown in equation (A.1), we have [using the notation $\text{o}(n^{-1}, N^{-1})$ to state that a quantity is an o of both $n$ and $N$, independently of how they mutually go to $\infty$]

$$r_N = 1 + \frac{1}{48}\left[\frac{1}{n} + \frac{1}{N}\right] + \text{o}\left(\frac{1}{n}, \frac{1}{N}\right).$$

Once again, it is not a troublesome quantity—given our objective of a 2/3 rate—since $|r_N - 1|$ goes to 0 like $1/N$. Nevertheless, it will sometimes simplify our work to keep it in the analyses.

Of central interest will be $(\kappa_N/\sigma_{n,N}^3)^{1/6}\hat{f}^{-1/4}$. We show in equation (A.3), that, if $\xi = \xi_2 + \varepsilon_N$, $\varepsilon_N$ small with respect to 1, there exists a sequence tending to a (finite) limit, denoted $\eta_N$, such that

$$(\kappa_N/\sigma_{n,N}^3)^{1/6}\hat{f}^{-1/4}(\xi) = 1 - \tfrac{2}{5}\eta_N\varepsilon_N + \text{O}(\varepsilon_N^2).$$

Since $\varepsilon_N$ will be as big as $\text{O}(N^{-1/2})$, we will have to be more careful about this part in the final steps of the analysis.

Finally, the most problematic part will turn out to be $\text{Ai}(\kappa_N^{2/3}\zeta)$. We recall that it was shown in [11] that $\kappa_N^{2/3}\zeta \to s$ for $x_N(s) = \mu_{n,N} + \sigma_{n,N}s$. Here we will need to go one step further in the asymptotic expansion of $\kappa_N^{2/3}\zeta$. For ease of exposition, we will first focus on $x_N(s) = \mu_{n,N} + \sigma_{n,N}s$, and will get to $\tilde{x}_N(s) = \tilde{\mu}_{n,N} + \tilde{\sigma}_{n,N}s$—the "optimal" centering and scaling for $\varphi_\tau + \psi_\tau$—only after we understand what goes "wrong" in terms of rates with $x_N(s)$.

We will focus in the next several subsections on

$$\theta_{n,N}(x) = F_{n,N}(x)\left(\frac{\mu_{n,N}}{x}\right),$$

since the proof of Lemma 3, via Fact 2.2.1, will rest on our ability to analyze quantities of the type

$$\Delta_{n,N}(x_N(s)) \triangleq |\alpha_{n,N}\theta_{n,N}(x_N(s)) - \text{Ai}(s)|.$$

We split $\Delta_{n,N}$ into two parts:

$$\Delta_{n,N}(x_N(s)) \leq \alpha_{n,N}|\theta_{n,N}(x_N(s)) - r_N \text{Ai}(s)| + |r_N \alpha_{n,N} - 1||\text{Ai}(s)|$$
$$\leq \alpha_{n,N}\Delta^{\text{I}}_{n,N}(x_N(s)) + \Delta^{\text{II}}_{n,N}(s)$$

with

$$\Delta^{\text{I}}_{n,N}(x_N(s)) = |\theta_{n,N}(x_N(s)) - r_N \text{Ai}(s)|$$

and

$$\Delta^{\text{II}}_{n,N}(s) = |r_N \alpha_{n,N} - 1||\text{Ai}(s)|.$$

We will assume from now on that $\alpha_{n,N} = 1 + \text{O}(1/N)$. Of course, we will verify that the sequence we eventually use has this property. We are going to show that, for $s \in [s_0, \infty)$, and $C$ a function that might change from display to display but will always be continuous and nonincreasing,

(I1) $$N^{2/3}\Delta^{\text{I}}_{n,N}(x_N(s)) \leq C(s_0)\exp(-s/2),$$

(I2) $$N\Delta^{\text{II}}_{n,N}(s) \leq C(s_0)\exp(-s/2).$$

3.2.2. *Rationale for the rates.* Before delving into the details of finding the bounds, we want to explain what is the rationale behind the rates corresponding to the different elements of the upper bounding sum. Recall that we are now working with a particular centering and scaling, namely, $x_N(s) = \mu_{n,N} + \sigma_{n,N}s$:

- $\Delta^{\text{I}}_{n,N}$: essentially what happens is that, with this centering and scaling, $\kappa_N^{2/3}\zeta_N(x_N(s))$ converges to $s$ at speed $N^{2/3}$, and we are able to deal with the rest of the elements at this "speed."
- $\Delta^{\text{II}}_{n,N}$: here, of course, only the rate of convergence of $r_N \alpha_{n,N}$ to 1 matters. The proof of (I2) is an immediate consequence of the estimate we gave for $r_N$ and of our assumption that $\alpha_{n,N} = 1 + \text{O}(1/N)$, so we do not need to worry about it.

3.3. *Proof of* (I1). Given the dynamics of $\kappa_N^{2/3}\zeta(x_N(s))$ (explored in more detail in Appendix A.2.3), we can achieve our objective by using the coarse bound

$$N^{2/3}\Delta^{\text{I}}_{n,N}(x_N(s)) \leq N^{2/3}|\theta_{n,N}(x_N(s)) - r_N \text{Ai}(\kappa_N^{2/3}\zeta(x_N(s)))|$$
$$+ N^{2/3}r_N|\text{Ai}(\kappa_N^{2/3}\zeta(x_N(s))) - \text{Ai}(s)|.$$

We will use the following notation:

$$B_{n,N}(x_N(s)) = |\theta_{n,N}(x_N(s)) - r_N \text{Ai}(\kappa_N^{2/3}\zeta(x_N(s)))|,$$
$$D_{n,N}(x_N(s)) = \text{Ai}(\kappa_N^{2/3}\zeta(x_N(s))) - \text{Ai}(s).$$

Hence, the previous inequality is just

$$N^{2/3}\Delta_{n,N}^{\mathrm{I}} \leq N^{2/3}B_{n,N} + r_N N^{2/3}|D_{n,N}|,$$

where we have dropped the argument $x_N(s)$ in the interest of clarity. We will continue to do so when there is no ambiguity about the argument we are manipulating.

It will turn out that $B_{n,N}$ is more "robust" to recentering and rescaling than $D_{n,N}$: when modifying slightly the centering and scaling [i.e., going from $x_N(s)$ to $\tilde{x}_N(s)$], we will be able to achieve the same rate—$2/3$—for $B_{n,N}(\tilde{x}_N(s))$ but not for $D_{n,N}(\tilde{x}_N(s))$.

We are going to split $[s_0, \infty)$ into three varying intervals, namely, $[s_0, s_1]$, $[s_1, N^{1/6}s_1]$ and $[N^{1/6}s_1, \infty)$; on the first two of these intervals, we will principally use our understanding of $\zeta(\xi)$. This splitting will turn out to be natural because we will need to be precise on $[s_0, N^{1/6}s_1]$, and will rely on Taylor expansions to carry out the work. $N^{1/6}s_1$ is small enough at the relevant scale that they will be uniformly valid on $[s_0, N^{1/6}s_1]$. On the other hand, given the speed of decay of the Airy function, we will use coarse decay bounds for this special function on the last interval, where $s$ is necessarily large, because $N$ is large enough.

Before we delve into the details, let us jump ahead and explain what controls the speed of decay to zero of $B_{n,N}$ and $D_{n,N}$.

For $B_{n,N}(x_N(s))$, the key quantity is going to be

$$\left|\left(\frac{\kappa_N}{\sigma_{n,N}^3}\right)^{1/6}\hat{f}^{-1/4}(\xi)(x_N(s)/\mu_{n,N})^{-1} - 1\right|,$$

which, as we will see, is of "order" $\varepsilon_N(s)(=\xi(s) - \xi_2 = x_N(s)/\kappa_N - \xi_2)$. So $B_{n,N}$ is also of this order, and the announced rate will hold because, roughly, $\kappa_N^{2/3}\varepsilon_N(s)/s = \mathrm{O}(1)$. This is more a take-away, heuristic message than a precise mathematical statement, but with this in mind, we can make everything rigorous.

For $D_{n,N}(x_N(s))$, we will see that what really matters is the rate of convergence of

$$|\kappa_N^{2/3}\zeta(x_N(s)) - s|$$

on intervals where $\varepsilon_N(s)$ is "under control." The decay to zero of this quantity is what really hurts us in terms of rate when we cannot use (for the trade-off reasons described in Section 2.2) an optimal centering and scaling sequence for $D_{n,N}$, that is, when we have to work with $D_{n,N}(\tilde{x}_N(s))$ instead of $D_{n,N}(x_N(s))$. We will see in Section 3.5 how we can nonetheless overcome this difficulty.

3.3.1. *Bounds for* $N^{2/3}B_{n,N}$. Recall that $\theta_{n,N}(z) = F_{n,N}(z)(\mu_{n,N}/z)$. So, using equation (5), we get

$$B_{n,N}(x_N(s))$$

$$= |F_{n,N}(x_N(s))(\mu_{n,N}/x_N(s)) - r_N \mathrm{Ai}(\kappa_N^{2/3}\zeta)|$$

$$\leq r_N \left| \left( \frac{\kappa_N}{\sigma_{n,N}^3} \right)^{1/6} \hat{f}^{-1/4}(\xi)(\mu_{n,N}/x_N(s)) - 1 \right| |\mathrm{Ai}(\kappa_N^{2/3}\zeta)|$$

$$+ \mathrm{O}\left( \frac{1}{\kappa_N} \right) \left| \mathbf{E}^{-1}(\kappa_N^{2/3}\zeta)\mathbf{M}(\kappa_N^{2/3}\zeta) \left( \frac{\mu_{n,N}}{x_N(s)} \right) r_N \left( \frac{\kappa_N}{\sigma_{n,N}^3} \right)^{1/6} \hat{f}^{-1/4}(\xi) \right|.$$

Note that the last line [without $\mathrm{O}(\frac{1}{\kappa_N})$] was bounded by $C(s_0)e^{-s}$ on $[s_0, \infty)$ in [11], (A.8), so we just have to concentrate on the first part of the sum.

To do this, we split $[s_0, \infty)$ into $[s_0, s_1]$, $[s_1, N^{1/6}s_1]$ and $[N^{1/6}s_1, \infty)$. We treat the problems in decreasing order of difficulty.

*Case* $s \in I_{1,N} = [s_1, N^{1/6}s_1]$. First a note on $s_1$: it is chosen as in [11], (A.8), that is, it is such that, for $s \geq s_1$, $2/3\kappa_N^{2/3}\zeta \geq s$. As we explain in Appendix A.6.4, we can choose the same $s_1$ for all values of $\gamma$, the limit of $n/N$. Also, we can assume that $s_1 \geq 1$, which guarantees that Ai is positive on $[s_1, \infty)$, that Ai$'$ is negative and increasing on this interval, and that

$$\mathrm{Ai}(x) \leq \frac{\exp(-2x^{3/2}/3)}{2\pi^{1/2}},$$

using properties of the Airy function cited in [15], pages 393–394.

Now recall that we chose $x_N(s) = \mu_{n,N} + \sigma_{n,N}s$; in equation (A.3), it means that $\varepsilon_N(s) = s\sigma_{n,N}/\kappa_N$, since $\kappa_N\xi_2 = \mu_{n,N}$. As an aside, let us point out that in what follows, we will often drop the dependence on $s$ of $\varepsilon_N(s)$ to simplify the notation. We will also use the notation $\eta_N$ for a specified sequence of real numbers having a finite limit (when $\gamma$ is given) and defined in Appendix A.2. With this in mind, we have, through equation (A.3),

$$\left( \frac{\kappa_N}{\sigma_{n,N}^3} \right)^{1/6} \hat{f}^{-1/4}(\xi) \frac{\mu_{n,N}}{x_N(s)}$$

$$= \left( 1 - \frac{2}{5}\eta_N \varepsilon_N + \mathrm{O}(\varepsilon_N^2) \right) \left( 1 - \frac{\sigma_{n,N}s}{\mu_{n,N}} + \mathrm{O}(\varepsilon_N^2) \right)$$

$$= (1 + \tilde{\eta}_N \varepsilon_N + \mathrm{O}(\varepsilon_N^2)),$$

as $\kappa_N/\mu_{n,N} \asymp 1$. Once again, $\tilde{\eta}_N$ has a finite limit as $N \to \infty$. Now on $I_{1,N}$, $\varepsilon_N = \mathrm{O}(N^{-1/2})$. Note that this estimate is also valid if we use $x_N(s) = \tilde{\mu}_{n,N} + \tilde{\sigma}_{n,N}s$, under certain conditions on $\tilde{\mu}_{n,N}$ and $\tilde{\sigma}_{n,N}$. The ones we will eventually use easily satisfy those. Also, in what follows we use "coarse" O, that is, what appears in them may not be the optimal quantity, but it is definitely good enough for what we need. We conclude that

$$\forall s \in I_{1,N} \qquad \left(\frac{\kappa_N}{\sigma_{n,N}^3}\right)^{1/6} \hat{f}^{-1/4}(\xi)(\mu_{n,N}/x_N(s)) = 1 + \tilde{\eta}_N \varepsilon_N + \mathrm{O}(N^{-4/5}).$$

Hence, there exists $\chi(\gamma)$, a function independent of $N$ and $n$ and, hence, a constant at $\gamma$ fixed such that

$$N^{2/3} \left|\left(\frac{\kappa_N}{\sigma_{n,N}^3}\right)^{1/6} \hat{f}^{-1/4}(\xi)(\mu_{n,N}/x_N(s)) - 1\right| \leq \chi(\gamma) s$$

when $s \in I_{1,N}$, and $N$ is large enough. As explained in Appendix A.6, $\chi(\gamma)$ can be chosen independently of $\gamma$. Hence, we will call it simply $\chi$, which will represent a generic constant, which can change from display to display. [Note that $N(s_0, \gamma)$, the integer after which these estimates are valid, might depend on $\gamma$.]

Now since $s \geq s_1$, $2/3 \kappa_N \zeta^{3/2} \geq s$, and therefore,

$$\text{on } I_{1,N}, \qquad |\mathrm{Ai}(\kappa_N^{2/3} \zeta)| \leq \chi \exp(-s).$$

Combining the two previous results, we get, for $s \in I_{1,N}$,

$$N^{2/3} \left|\left(\frac{\kappa_N}{\sigma_{n,N}^3}\right)^{1/6} \hat{f}^{-1/4}(\xi)(\mu_{n,N}/x_N(s)) - 1\right| |\mathrm{Ai}(\kappa_N^{2/3} \zeta)| \leq \chi s \exp(-s)$$

$$\leq \chi \exp(-s/2).$$

*Case $s \in I_{2,N} = [N^{1/6} s_1, \infty)$.* Here we can act heavy-handedly: the fast decay of the Airy function alone will suffice for our purposes. As a matter of fact, we use the very coarse bound

$$\left|\left(\frac{\kappa_N}{\sigma_{n,N}^3}\right)^{1/6} \hat{f}^{-1/4}(\xi)(\mu_{n,N}/x_N(s)) - 1\right| |\mathrm{Ai}(\kappa_N^{2/3} \zeta)|$$

$$\leq \left|\left(\frac{\kappa_N}{\sigma_{n,N}^3}\right)^{1/6} \hat{f}^{-1/4}(\xi)(\mu_{n,N}/x_N(s)) \mathrm{Ai}(\kappa_N^{2/3} \zeta)\right| + |\mathrm{Ai}(\kappa_N^{2/3} \zeta)|.$$

Now, since $\mathrm{Ai} \leq \mathbf{M} \mathbf{E}^{-1}$, using (A.8) in [11], it follows that

$$\left|\left(\frac{\kappa_N}{\sigma_{n,N}^3}\right)^{1/6} \hat{f}^{-1/4}(\xi)(\mu_{n,N}/x_N(s)) \mathrm{Ai}(\kappa_N^{2/3} \zeta)\right| \leq \chi \exp(-s).$$

Since $s \geq N^{1/6} s_1$, we have $s \geq 4/3 \log(N)$, which implies that

$$N^{2/3} \exp(-s) \leq \exp(-s/2).$$

Similarly (using $2/3\kappa_N\zeta^{3/2} \geq s$, and the bound on $\mathrm{Ai}(x)$ given in [15], page 394), $|\mathrm{Ai}(\kappa_N^{2/3}\zeta)| \leq \exp(-2\kappa_N\zeta^{3/2}/3) \leq e^{-s} \leq N^{-2/3}e^{-s/2}$. We can therefore conclude that

$$N^{2/3}r_N\left|\left(\frac{\kappa_N}{\sigma_{n,N}^3}\right)^{1/6}\hat{f}^{-1/4}(\xi)(\mu_{n,N}/x_N(s)) - 1\right||\mathrm{Ai}(\kappa_N^{2/3}\zeta)| \leq \chi(e^{-s/2} + e^{-s})$$

$$\leq \chi\exp(-s/2).$$

*Case* $s \in I_{0,N} = [s_0, s_1]$. The arguments we used on $I_{1,N}$ apply and show that

$$N^{2/3}r_N\left|\left(\frac{\kappa_N}{\sigma_{n,N}^3}\right)^{1/6}\hat{f}^{-1/4}(\xi)(x_N(s)/\mu_{n,N})^{-1} - 1\right| \leq \chi(|s| \vee 1).$$

For $s \in I_{0,N}$, the sequence $\kappa_N^{2/3}\zeta(x_N(s))$ is bounded, because, at $n, N$ fixed, $\zeta$ is an increasing function of $s$, and it is uniformly bounded in $N$ since $\kappa_N^{2/3}\zeta(x_N(s_0)) \to s_0$ and $\kappa_N^{2/3}\zeta(x_N(s_1)) \to s_1$. So we have, for $N$ large enough,

$$\forall s \in I_{0,N} \quad -2(|s_0| \vee 1) \leq \kappa_N^{2/3}\zeta(x_N(s)) \leq 2s_1,$$

and therefore,

$$\forall s \in I_{0,N} \quad |\mathrm{Ai}(\kappa_N^{2/3}\zeta)| \leq \sup_{s \in [-2(|s_0|\vee 1), 2s_1]}|\mathrm{Ai}(s)|.$$

A point on notation: in several displays, we will use quantities like $-2|s_0|$. Our arguments are correct if $s_0 \neq 0$, but we can replace everywhere $|s_0|$ by $|s_0| \vee 1$ and they are valid whether or not $s_0 = 0$. Aware of this, we nevertheless choose to write $|s_0|$ instead of $|s_0| \vee 1$ to avoid cumbersome notation, and given the fact that this is really a minor detail. Going back to the analysis, we obtain

$$N^{2/3}r_N\left|\left(\frac{\kappa_N}{\sigma_{n,N}^3}\right)^{1/6}\hat{f}^{-1/4}(\xi)(x_N(s)/\mu_{n,N})^{-1} - 1\right||\mathrm{Ai}(\kappa_N^{2/3}\zeta)|$$

$$\leq \chi(|s| \vee 1)\sup_{s \in [-2|s_0|, 2s_1]}|\mathrm{Ai}(s)| \leq C(s_0)\exp(-s/2).$$

We can choose for $C$ in the previous display, up to a constant, the function $C(s_0) = \{\sup_{s\in[-2|s_0|,2s_1]}|\mathrm{Ai}(s)|\sup_{s\in[-2|s_0|,2s_1]}[\exp(s/2)(|s| \vee 1)]\}$. This is of course a continuous and nonincreasing function of $s_0$. $C$ (and the aforementioned constant) can be chosen independently of $\gamma$ through arguments similar to the ones mentioned when the problem of the dependence of $\chi$ on $\gamma$ arose on $I_{1,N}$.

We have therefore shown that if $n/N \to \gamma \in [1, \infty)$, for all $s_0$ we can find $N(s_0, \gamma)$ such that for $s \geq s_0$ and $N > N(s_0, \gamma)$,

$$N^{2/3} B_{n,N}(x_N(s)) \leq C(s_0) e^{-s/2}.$$

$C$ is a continuous nonincreasing function that is constant on $[s_1, \infty)$ and independent of $\gamma$. [Actually, the explicit (up to a constant) expression for $C$ given above shows that we could replace $C$ by a constant, but it is of little use given what follows.]

3.3.2. *Bounding* $N^{2/3}|D_{n,N}(x_N(s))|$. Here we are of course going to be relying heavily on first-order asymptotics for $\kappa_N^{2/3}\zeta(x_N(s))$, given in equation (A.4) and the fact that $|\text{Ai}'(s)|$ is decreasing on $\mathbb{R}_+$ (Ai' is increasing because Ai'' is positive on $\mathbb{R}_+$ according to the Airy equation, and we know that Ai' is negative on $\mathbb{R}_+$). We use the same decomposition of $[s_0, \infty)$ as in the $B_{n,N}$ case.

*Case* $s \in I_{1,N}$. We have from equation (A.4)

$$\kappa_N^{2/3}\zeta(x_N(s)) = \frac{\varepsilon_N \kappa_N}{\sigma_{n,N}}\left(1 + \frac{2}{5}\varepsilon_N \eta_N + O(\varepsilon_N^2)\right).$$

Our choice of centering and scaling gives $\varepsilon_N = \varepsilon_N(s) = s\sigma_{n,N}/\kappa_N$. So, when $N$ is large enough, $\min(\kappa_N^{2/3}\zeta(x_N(s)), s) \geq s/2$. Therefore, using the mean-value theorem and the fact that $|\text{Ai}'|$ is decreasing,

$$|\text{Ai}(\kappa_N^{2/3}\zeta) - \text{Ai}(s)| \leq |\text{Ai}'(s/2)||\kappa_N^{2/3}\zeta - s| \leq \chi(\gamma)|\text{Ai}'(s/2)|s^2 \frac{\sigma_{n,N}}{\kappa_N}.$$

Now since $s \geq s_1 \geq 1$, using the first formula on page 394 of [15], we get

$$|\text{Ai}'(s/2)|s^2 \leq \frac{s^{9/4} e^{-s^{3/2}/(3\sqrt{2})}}{\sqrt{\pi} 2^{1/4}} \leq \chi e^{-s/2},$$

which shows that indeed $N^{2/3}|D_{n,N}| \leq \chi e^{-s/2}$ on $I_{1,N}$. ($\chi$ can be chosen independently of $\gamma$, as shown by arguments similar to those made in Appendix A.6.)

*Case* $s \in I_{2,N}$. In this case, we can be as rough as we were on the corresponding interval for $B_{n,N}$:

$$N^{2/3}|D_{n,N}| \leq N^{2/3}(|\text{Ai}(\kappa_N^{2/3}\zeta)| + |\text{Ai}(s)|) \leq \chi e^{-s/2}.$$

*Case* $s \in I_{0,N}$. Here, as we showed in the corresponding case for $B_{n,N}$, $\kappa_N^{2/3}\zeta$ stays uniformly bounded when $N$ goes to $\infty$. Moreover, the arguments given for $I_{1,N}$ still hold, and we have

$$N^{2/3}|\text{Ai}(\kappa_N^{2/3}\zeta) - \text{Ai}(s)| \leq \chi s^2 \sup_{[-2|s_0|, 2s_1]} |\text{Ai}'(s)| \leq C(s_0) e^{-s/2}.$$

$C$ can be chosen—up to a constant—to be the continuous nonincreasing function $C(s_0) = \{\sup_{s\in[-2|s_0|,2s_1]} |\mathrm{Ai}'(s)| \sup_{s\in[-2|s_0|,2s_1]}[\exp(s/2)(s^2 \vee 1)]\}$. Note that this $C$ tends to $\infty$ as $s_0 \to -\infty$, and does it like $|s_0|^{1/4}$, using well-known properties of $\mathrm{Ai}'(s)$. The constants $\chi$ (which look like they are functions of $\gamma$) can be shown to be independent of $\gamma$, through estimates given in Appendix A.6.

We have shown the inequality (I1).

3.4. *About* (P2), (P3) *and* (P1). We conclude from the previous subsection that if we chose $\tilde{\mu}_{n,N} = \mu_{n-1,N}$ and $\tilde{\sigma}_{n,N} = \sigma_{n-1,N}$, we would have $N^{2/3}\|G_\tau - G\|_2 \leq C(s_0)\exp(-s/2)$. The problems for (P3) and (P1) are essentially the same, so we will just focus on (P3). We have to understand how a nonoptimal centering affects the rate of convergence of $H_\tau$ to $G$. What we will see in Section 3.4.2 is that using $\mu_{n-1,N}$ and $\sigma_{n-1,N}$ on $H_\tau$ (instead of the "optimal" $\mu_{n,N-1}$ and $\sigma_{n,N-1}$) prevents us from being able to show convergence of $H_\tau$ to $G$ at a speed faster than $N^{1/3}$. Since we will need to compromise between $H_\tau$ and $G_\tau$ (because what is optimal for one is not optimal for the other), "favoring" $G_\tau$ over $H_\tau$ (or vice-versa) turns out to be a bad choice for the global problem and we need to investigate a new problem we now set up and will solve in Section 3.5.

Let $\tilde{\mu}_{n,N}$ and $\tilde{\sigma}_{n,N}$ be, respectively, a centering and scaling sequence. Let us further assume that $\tilde{\mu}_{n,N} - \mu_{n-1,N} = \mathrm{O}(1)$ and $(\tilde{\sigma}_{n,N}/\sigma_{n-1,N} - 1) = \mathrm{O}(N^{-1})$. We will show that, in this situation, (P2) holds. Later, we will show that we can find a pair $(\tilde{\mu}_{n,N}, \tilde{\sigma}_{n,N})$ such that (P1) holds and at this point Theorem 2 will be proven.

So we are now dealing with $\tilde{x}_N(s) = \tilde{\mu}_{n,N} + \tilde{\sigma}_{n,N}s$. We have

$$\varepsilon_{n-1,N} = (\tilde{\mu}_{n,N} - \mu_{n-1,N})/\kappa_{n-1,N} + s\tilde{\sigma}_{n,N}/\kappa_{n-1,N}.$$

Hence, when we work with $\kappa_{n-1,N}\varepsilon_{n-1,N}/\sigma_{n-1,N}$, the first term is already $\mathrm{O}(N^{-1/3})$. This is fundamentally what "harms" the rate when we do not choose an "optimal" centering and scaling sequence.

We now look in more detail at $B_{n-1,N}(\tilde{x}_N(s))$ and $D_{n-1,N}(\tilde{x}_N(s))$. We will show in Appendix A.1.4, that $\alpha_{n-1,N} = \sqrt{a_N}\sigma_{n-1,N}^{1/2}\tilde{\sigma}_{n,N}/\mu_{n-1,N}$ goes to 1 at rate 1 [i.e., it is $1 + \mathrm{O}(N^{-1})$] and that it therefore satisfies the assumptions put forward in Section 3.2.2.

3.4.1. $B_{n-1,N}(\tilde{x}_N(s))$ *situation*. As we will soon see, this part is not really problematic: we keep the $2/3$ rate even when choosing a—reasonable—nonoptimal sequence.

Note that $I_{2,N}$ is not a problem, as the upper bounding relied on the speed of the decay of the Airy function, and using (A.4), we have $\kappa_{n-1,N}^{2/3}\zeta(\tilde{x}_N(s)) \geq 0.5 + 3s/4$, if $s \geq s_1$ and $n, N$ large enough.

Now on $I_{1,N}$ and $I_{0,N}$, we are still fine: using the fact that $\varepsilon_{n-1,N} = (\tilde{\mu}_{n,N} - \mu_{n-1,N})/\kappa_{n-1,N} + s\tilde{\sigma}_{n,N}/\kappa_{n-1,N}$ and $(\tilde{\mu}_{n,N} - \mu_{n-1,N})/\kappa_{n-1,N} = \mathrm{O}(N^{-1})$, the analysis carried above still holds (we give more details on this in Section 3.5), and we have, $\forall s \in I_{1,N}$,

$$r_{n-1,N} \left(\frac{\kappa_{n-1,N}}{\sigma_{n-1,N}^3}\right)^{1/6} \hat{f}^{-1/4}(\xi) \frac{\mu_{n-1,N}}{\tilde{x}_N(s)}$$
$$= 1 - \frac{2}{5}\eta_N^{\{n-1\}}\varepsilon_{n-1,N} - \frac{\tilde{\sigma}_{n,N}s}{\mu_{n-1,N}} + \mathrm{O}(N^{-4/5})$$
$$= 1 + \tilde{\eta}_N^{\{n-1\}}\varepsilon_{n-1,N} + \mathrm{O}(N^{-4/5}).$$

Using the same ideas as above, we conclude that

$$N^{2/3} B_{n-1,N}(\tilde{x}_N(s)) \leq C(s_0)\exp(-s/2),$$

where, once again, $C$ can be chosen to be continuous and nonincreasing.

3.4.2. *$D_{n-1,N}(\tilde{x}_N(s))$ situation.* The problem is essentially the same as it was in Section 3.3.2: dealing now with $\varepsilon_{n-1,N}$, we get, by simply applying (A.4)

$$\kappa_{n-1,N}^{2/3}\zeta(\tilde{x}_N(s)) = s + \mathrm{O}(N^{-1/3}).$$

We can then apply the mean-value theorem to get exponential bounds, but the problem remains: the speed cannot be shown to be faster than $N^{1/3}$.

This completes the proof of (P2), since it shows that (3) holds true. It also shows that (P3) holds.

3.5. *Better centering and scaling.* We now turn to the problem of finding centering and scaling sequences such that (P1) holds and the assumptions we made (about the relationship between the sequences $\tilde{\mu}_{n,N}$ and $\tilde{\sigma}_{n,N}$ on one hand, and $\mu_{n-1,N}$ and $\sigma_{n-1,N}$ on the other) in the previous subsections are valid. What we need to do is relatively clear: we need to compensate the $N^{-1/3}$ deviation of $G_\tau$ from $G$ by the $N^{-1/3}$ deviation of $H_\tau$ from $G$ to get a higher rate of convergence for $G_\tau + H_\tau$ to $2G$. The discussions above show that the only region that is problematic is $I_{0,N} \cup I_{1,N}$, and the problems arise only for $D_{n-1,N}$ and $D_{n,N-1}$.

So let us now focus on

$$\Delta_{n,N}^F(\tilde{x}_N(s)) = |\varphi_\tau(\tilde{x}_N(s)) + \psi_\tau(\tilde{x}_N(s)) - \sqrt{2}\mathrm{Ai}(s)|.$$

We have

$$\sqrt{2}\Delta_{n,N}^F(\tilde{x}_N(s))$$
(6)
$$\leq \alpha_{n-1,N}|B_{n-1,N}(\tilde{x}_N(s))| + \alpha_{n,N-1}|B_{n,N-1}(\tilde{x}_N(s))|$$
$$+ |\alpha_{n-1,N}r_{n-1,N}D_{n-1,N}(\tilde{x}_N(s)) + \alpha_{n,N-1}r_{n,N-1}D_{n,N-1}(\tilde{x}_N(s))|$$

and our aim is to show that (with the same quantifiers as usual)

(I3) $\quad N^{2/3}\Delta_{n,N}^F(\tilde{x}_N(s)) \leq C(s_0)\exp(-s/2) \quad$ on $[s_0,\infty)$.

Once again, $C$ will be shown to be a continuous and nonincreasing function.

Neither $B_{n-1,N}(\tilde{x}_N(s))$ nor $B_{n,N-1}(\tilde{x}_N(s))$ will cause problems in terms of rates (as shown by the arguments given in Section 3.4.1), so we just need to focus on the first term of (6). Before we proceed to giving the needed explanations, let us introduce the notation:

$$\tilde{\varepsilon}_N = \frac{\tilde{\mu}_{n,N} - \mu_{n-1,N}}{\kappa_{n-1,N}} + \frac{\tilde{\sigma}_{n,N}}{\kappa_{n-1,N}}s,$$

$$\tilde{\varepsilon}_{N-1} = \frac{\tilde{\mu}_{n,N} - \mu_{n,N-1}}{\kappa_{n,N-1}} + \frac{\tilde{\sigma}_{n,N}}{\kappa_{n,N-1}}s,$$

$$\delta_{n-1,N} \triangleq \frac{\tilde{\mu}_{n,N} - \mu_{n-1,N}}{\kappa_{n-1,N}},$$

$$\delta_{n,N-1} \triangleq \frac{\tilde{\mu}_{n,N} - \mu_{n,N-1}}{\kappa_{n,N-1}}.$$

We have $\tilde{x}_N(s)/\kappa_{n-1,N} = \xi_2^{(n-1,N)} + \tilde{\varepsilon}_N$ and similarly for $\tilde{\varepsilon}_{N-1}$. From now on we make the assumptions that $\tilde{\mu}_{n,N} - \mu_{n-1,N} = \mathrm{O}(1)$, $\tilde{\mu}_{n,N} - \mu_{n,N-1} = \mathrm{O}(1)$ and that $\tilde{\sigma}_{n,N}/\sigma_{n-1,N} - 1, \tilde{\sigma}_{n,N}/\sigma_{n,N-1} - 1 = \mathrm{O}(N^{-1}, n^{-1})$, which we will show (in Appendices A.1.2 and A.1.3) hold for our eventual choice of centering and scaling. Finally, we note that since $r_{n-1,N}$ and $r_{n,N-1}$ are $1 + \mathrm{O}(1/N)$, they will not affect the discussion that follows, and so we drop them for the sake of simplicity.

The only real problem is with

$$\tilde{\Delta}_{n,N}^F(\tilde{x}_N(s)) \triangleq |\alpha_{n-1,N}D_{n-1,N}(\tilde{x}_N(s)) + \alpha_{n,N-1}D_{n,N-1}(\tilde{x}_N(s))|.$$

3.5.1. *Analysis of* $\tilde{\Delta}_{n,N}^F$. Let us first remark that the coarse approach explained in detail above showed that on $I_{2,N}$ we have $N^{2/3}D_{n-1,N}(\tilde{x}_N(s)) \leq C(s_0)e^{-s/2}$, $N^{2/3}D_{n,N-1}(\tilde{x}_N(s)) \leq C(s_0)e^{-s/2}$ and, therefore, we do not have to worry about this part of the real line. In what follows, we therefore assume that $s \in I_{0,N} \cup I_{1,N}$.

To show that $\tilde{\Delta}_{n,N}^F$ goes to zero at the $2/3$ rate claimed in (I3), we just have to use Taylor's formula with integral remainder: we are going to choose $\tilde{\mu}_{n,N}$ such that the first-order terms [in $\mathrm{Ai}'(s)$] that appear cancel out. Then, we will show that the remainder is $\mathrm{O}(N^{-2/3}e^{-s/2})$.

Using equation (A.4), we have

$$\kappa_{n-1,N}^{2/3}\zeta(\tilde{x}_N(s)) = \frac{\tilde{\varepsilon}_N\kappa_{n-1,N}}{\sigma_{n-1,N}}\left(1 + \frac{2}{5}\tilde{\varepsilon}_N\eta_N + \mathrm{O}(\tilde{\varepsilon}_N^2)\right).$$

Let us call
$$u_{n-1,N}(s) \triangleq \kappa_{n-1,N}^{2/3}\zeta(\tilde{x}_N(s)) - s.$$

While keeping in mind that $u_{n-1,N}(s)$ depends on $s$, we will often drop the $s$ to alleviate the notation. Now we use the fact that $\text{Ai}''(x) = x\text{Ai}(x)$ to get (through Taylor's formula with integral remainder)

$$\text{Ai}(\kappa_{n-1,N}^{2/3}\zeta(\tilde{x}_N(s))) = \text{Ai}(s) + u_{n-1,N}\text{Ai}'(s)$$
$$+ \int_0^{u_{n-1,N}} (u_{n-1,N} - y)(s+y)\text{Ai}(s+y)\,dy,$$

and therefore,

$$\alpha_{n-1,N}D_{n-1,N} + \alpha_{n,N-1}D_{n,N-1}$$
$$= (\alpha_{n-1,N}u_{n-1,N} + \alpha_{n,N-1}u_{n,N-1})\text{Ai}'(s)$$
$$+ \alpha_{n-1,N}R_N(s) + \alpha_{n,N-1}R_{N-1}(s).$$

We start by focusing on $(\alpha_{n-1,N}u_{n-1,N} + \alpha_{n,N-1}u_{n,N-1})$ and we will show later that $R_N$ and $R_{N-1}$ decay to zero fast enough for our needs. Since $\alpha_{n-1,N} = 1 + \text{O}(1/N)$, it is clear that if we can show that $R_N$ and $R_{N-1}$ decay to zero at rate $2/3$, $\alpha_{n-1,N}R_N$ and $\alpha_{n,N-1}R_{N-1}$ will decay to zero at the same speed.

Before we proceed, we recall that $\alpha_{n-1,N} = \sqrt{a_N}\sigma_{n-1,N}^{1/2}\tilde{\sigma}_{n,N}/\tilde{\mu}_{n,N}$.

*Remark on $u_{n-1,N}$.* Using equation (A.4), we have
$$u_{n-1,N}(s) = \frac{\kappa_{n-1,N}}{\sigma_{n-1,N}}\tilde{\varepsilon}_N\left(1 + \frac{2}{5}\tilde{\varepsilon}_N\eta_N + \text{O}(\tilde{\varepsilon}_N^2)\right) - s,$$

where $\eta_N$ has a finite limit. Recall that $\tilde{\varepsilon}_N = \delta_{n-1,N} + s\tilde{\sigma}_{n,N}/\kappa_{n-1,N}$, where $\delta_{n-1,N} = \text{O}(N^{-1})$. So

$$u_{n-1,N}(s) = \frac{\kappa_{n-1,N}}{\sigma_{n-1,N}}\delta_{n-1,N} + \left(\frac{\tilde{\sigma}_{n,N}}{\sigma_{n-1,N}} - 1\right)s$$
$$+ \frac{2\eta_N\kappa_{n-1,N}}{5\sigma_{n-1,N}}\tilde{\varepsilon}_N^2 + \text{O}\left(\frac{\tilde{\varepsilon}_N^3\kappa_{n-1,N}}{\sigma_{n-1,N}}\right).$$

Let us focus on $\tilde{\varepsilon}_N$ for a moment. It is clear that $|\tilde{\varepsilon}_N| \leq \chi N^{-2/3}(|s| \vee 1)$, where $\chi$ is independent of $N$, $s$ and $\gamma$. It follows that

$$N^{2/3}\left|\frac{2\eta_N\kappa_{n-1,N}}{5\sigma_{n-1,N}}\tilde{\varepsilon}_N^2 + \text{O}(\tilde{\varepsilon}_N^3\kappa_{n-1,N}/\sigma_{n-1,N})\right| \leq \chi(|s|^3 \vee 1).$$

The fact that $\chi$ can be chosen independently of $\gamma$ is shown in Appendix A.6. In other respects, we recall that, as mentioned in [15] for $x > 0$,

$$|\text{Ai}'(x)| \leq \frac{x^{1/4}e^{-2/3x^{3/2}}}{2\pi^{1/2}}\left(1 + \frac{7}{48x^{3/2}}\right),$$

from which we deduce, along the lines of the analyses we did before, that

$$N^{2/3}\left|\frac{2\eta_N\kappa_{n-1,N}}{5\sigma_{n-1,N}}\tilde{\varepsilon}_N^2 + \mathrm{O}(\tilde{\varepsilon}_N^3\kappa_{n-1,N}/\sigma_{n-1,N})\right||\mathrm{Ai}'(s)| \leq C(s_0)e^{-s/2},$$

$C(s_0) = \{\sup_{[s_1,\infty)}[e^{s/2-2s^{3/2}/3}s^{1/4}] + \sup_{[-2|s_0|,2s_1]} e^{s/2}|\mathrm{Ai}'(s)|(|s|^3 \vee 1)\}$, for instance. Hence, $C$ can be chosen to be a continuous and nonincreasing function; in this particular instance we could replace it by a constant.

In other words, for our rates purposes, it is enough to focus on

$$\tilde{u}_{n-1,N}(s) = \frac{\kappa_{n-1,N}}{\sigma_{n-1,N}}\delta_{n-1,N} + \left(\frac{\tilde{\sigma}_{n,N}}{\sigma_{n-1,N}} - 1\right)s = \mathrm{O}((|s| \vee 1)N^{-1/3}).$$

Note that the same analysis applies to $u_{n,N-1}$.

*An intuitive choice for $\tilde{\mu}_{n,N}$.* At an intuitive level, our biggest problem comes from the "centering" problem, so it is natural to try to get rid of it by canceling its effect and then verify that we then get the rate we were expecting.

The "centering" term, the one that appears because one time (i.e., for $\varphi$-related matters, or parameters $(n-1, N)$] the "optimal" centering is $\mu_{n-1,N}$ and the other time [i.e., for $\psi$-related matters, or parameters $(n, N-1)$] it is $\mu_{n,N-1}$, is

$$\begin{aligned}
\mathbf{c_N} &\triangleq \alpha_{n-1,N}\frac{\kappa_{n-1,N}}{\sigma_{n-1,N}}\delta_{n-1,N} + \alpha_{n,N-1}\frac{\kappa_{n,N-1}}{\sigma_{n,N-1}}\delta_{n,N-1} \\
&= \tilde{\sigma}_{n,N}\sqrt{a_N}\left[\left(\frac{\tilde{\mu}_{n,N}}{\mu_{n-1,N}} - 1\right)\sigma_{n-1,N}^{-1/2} + \left(\frac{\tilde{\mu}_{n,N}}{\mu_{n,N-1}} - 1\right)\sigma_{n,N-1}^{-1/2}\right].
\end{aligned}$$
(7)

Canceling it leads us to choose

(centering)
$$\tilde{\mu}_{n,N} = \left(\frac{1}{\sigma_{n-1,N}^{1/2}} + \frac{1}{\sigma_{n,N-1}^{1/2}}\right)$$
$$\times \left(\frac{1}{\mu_{n-1,N}\sigma_{n-1,N}^{1/2}} + \frac{1}{\mu_{n,N-1}\sigma_{n,N-1}^{1/2}}\right)^{-1}.$$

We will show in Appendix A.1.2 that we indeed have $\tilde{\mu}_{n,N} - \mu_{n-1,N} = \mathrm{O}(1)$.

*Study of $\alpha_{n-1,N}u_{n-1,N} + \alpha_{n,N-1}u_{n,N-1}$.* The conclusion of our remark on $u_{n-1,N}$ was that we can focus on $\tilde{u}_{n-1,N}$ rather than $u_{n-1,N}$ when dealing with the 2/3 rate. We have already seen that the choice made in (centering) led to

(8)
$$\begin{aligned}
\mathbf{s_N}s &\triangleq \alpha_{n-1,N}\tilde{u}_{n-1,N}(s) + \alpha_{n,N-1}\tilde{u}_{n,N-1}(s) - c_N \\
&= \sqrt{a_N}\tilde{\sigma}_{n,N}\left\{\frac{\sigma_{n-1,N}^{1/2}}{\mu_{n-1,N}}\left(\frac{\tilde{\sigma}_{n,N}}{\sigma_{n-1,N}} - 1\right) + \frac{\sigma_{n,N-1}^{1/2}}{\mu_{n,N-1}}\left(\frac{\tilde{\sigma}_{n,N}}{\sigma_{n,N-1}} - 1\right)\right\}s.
\end{aligned}$$

So a reasonable choice is to pick $\tilde{\sigma}_{n,N}$ so as to cancel this term, that is, after defining

$$\gamma_{n,N} = \frac{\mu_{n-1,N}\sigma_{n,N-1}^{1/2}}{\mu_{n,N-1}\sigma_{n-1,N}^{1/2}},$$

(scaling) $$\tilde{\sigma}_{n,N} = (1+\gamma_{n,N})\left(\frac{1}{\sigma_{n-1,N}} + \frac{\gamma_{n,N}}{\sigma_{n,N-1}}\right)^{-1}.$$

*Remark on the centering and scaling.* Note also that by picking a $\tilde{\sigma}_{n,N}$ such that the "scaling sequence" $s_N$ is an $\mathrm{O}(N^{-2/3})$, and a $\tilde{\mu}_{n,N}$ that makes the "centering sequence" $c_N$ be $\mathrm{O}(N^{-2/3})$, we get that

$$(\alpha_{n-1,N}u_{n-1,N}(s) + \alpha_{n,N-1}u_{n,N-1}(s)) = (|s| \vee 1)\mathrm{O}(N^{-2/3})$$

and we obtain the $N^{2/3}$ speed for the original problem.

*Bounding the remainders appearing in the Taylor expansion.* We first need to remark that

$$\text{On } I_{1,N}, \qquad u_{n-1,N}(s) = \mathrm{O}(s^3 N^{-1/3}),$$

$$\text{On } I_{0,N}, \qquad u_{n-1,N}(s) = \mathrm{O}(N^{-1/3})_{s_0}.$$

We use the notation $\mathrm{O}(\cdot)_{s_0}$ to emphasize the fact that the constant implicit in the O possibly depends on $s_0$. These estimates also apply to $u_{n,N-1}$. All the implicit constants can be chosen uniformly with respect to $\gamma$.

On $I_{1,N}$, we have, if $N$ is large enough, $s + u_{n-1,N}(s) \geq s/2$. Hence, we get, using the fact that the Airy function is nonincreasing and positive on this interval,

$$R_N(s) = \left|\int_0^{u_{n-1,N}} (u_{n-1,N} - y)(s+y)\mathrm{Ai}(s+y)\,dy\right|$$

$$\leq (|u_{n-1,N}(s)|^2/2)((s + u_{n-1,N}(s)) \vee s)\mathrm{Ai}(s/2).$$

Now $((s + u_{n-1,N}(s)) \vee s) = \mathrm{O}(s^3)$, and we therefore get, using only the control on $N$ provided by $u_{n-1,N}(s)^2$,

$$N^{2/3}R_N(s) \leq \chi\exp(-s/2).$$

On the other hand, on $I_{0,N}$, we have, for $N$ large enough, $-2|s_0| \leq s + u_{n-1,N}(s) \leq 2s_1$, and so we bound the remainder by

$$R_N(s) = \left|\int_0^{u_{n-1,N}} (u_{n-1,N} - y)(s+y)\mathrm{Ai}(s+y)\,dy\right|$$

$$\leq (|u_{n-1,N}(s)|^2/2)2(s_1 \vee |s_0|)\max_{[-2|s_0|,2s_1]}|\mathrm{Ai}(s)|.$$

Again, using the fact that $u_{n-1,N}(s) = \mathrm{O}(N^{-1/3})_{s_0}$ on this interval, we conclude

$$N^{2/3} R_N(s) \leq C(s_0) \exp(-s/2).$$

*Conclusion.* The same analysis applies to $R_{N-1}$, and this finishes the proof of the fact that

$$N^{2/3} \tilde{\Delta}_{n,N}^F(\tilde{x}_N(s)) \leq C(s_0) \exp(-s/2),$$

with a choice of centering and scaling satisfying the conditions set forth in our remark on centering and scaling.

## 4. Discussion.

4.1. *Centering and scaling.* The proof confirmed the empirical findings in [11] that small changes to $\tilde{\mu}_{n,N}$ and $\tilde{\sigma}_{n,N}$ can drastically improve the quality of the Tracy–Widom approximation, and the relevance of the Tracy–Widom law in small samples. This is an important fact for applications.

We recall the main conclusion of the analysis that we carried above: there is some liberty in choosing the centering and the scaling, as long as the chosen centering and scaling sequences (resp. $\tilde{\mu}_{n,N}$ and $\tilde{\sigma}_{n,N}$) satisfy the following properties [see equations (7) and (8) for the definitions of $c_N$ and $s_N$]:

$$c_N = \mathrm{O}(N^{-2/3}),$$
$$s_N = \mathrm{O}(N^{-2/3})$$

and $\tilde{\mu}_{n,N} - \mu_{n-1,N} = \mathrm{O}(1)$, $\tilde{\sigma}_{n,N}/\sigma_{n-1,N} = 1 + \mathrm{O}(N^{-1})$.

We finally note that different choices of $\tilde{\mu}_{n,N}$ and $\tilde{\sigma}_{n,N}$ (different from the ones indicated in Theorem 2, but satisfying the conditions just mentioned) might affect how the bounding functions behave with respect to $\gamma$. In other words, we might not be guaranteed that the bounding function $M$ in Theorem 2 can be chosen to be independent of $\gamma$, but we are guaranteed that, for any chosen $\gamma$, the convergence is happening at rate $2/3$. We give more detail on this issue in the next subsection and want to point out that with our choice of $\tilde{\mu}_{n,N}$ and $\tilde{\sigma}_{n,N}$ the bounding function $M$ in Theorem 2 can be chosen to be independent of $\gamma$.

4.2. *The bounding function $M$.* In Theorem 2, a bounding function $M$ appears; it is important to know how it behaves. We need a few preliminaries to ease the discussion.

4.2.1. *A remark on $\|\bar{S}\|_1(s)$.* First, we need to understand the behavior of $\|\bar{S}\|_1$ a as function of $s$, when $\bar{S}$ is viewed as an operator on $L^2([s,\infty))$.

Using the integral representation of $\bar{S}$, we see that it is a positive symmetric operator on $L^2([s,\infty))$, and so we can evaluate, through Mercer's theorem, its trace class norm (see, e.g., [13], Chapter 30 and, in particular, Section 30.5, which can be easily extended to intervals of the type $[a,\infty)$ when the kernels have exponential decay properties at $\infty$ similar to that of the Airy operator). Using Airy integrals spelled out in [14], AI.11(iv), we have for $\bar{S}$, as an operator on $[s,\infty)$,

$$\|\bar{S}\|_1(s) = \int_0^\infty \int_s^\infty \mathrm{Ai}^2(x+u)\,dx\,du$$
$$= (-\mathrm{Ai}(s)\mathrm{Ai}'(s) - 2s\mathrm{Ai}'(s)^2 + 2s^2\mathrm{Ai}^2(s))/3.$$

We can show that this quantity is unbounded as $s \to -\infty$ by using elementary properties of the Airy function. As an aside, let us mention that $\|\bar{S}\|_1(s)$ is of course a decreasing function of $s$. The order of magnitude of the previous quantity is $|s|^{3/2}$ when $s \to -\infty$. On the other hand, the speed of decay of the Airy function at $\infty$ shows that, on $[s_0,\infty)$,

$$\|\bar{S}\|_1(s) \leq C(s_0)\exp(-s),$$

for a continuous and nonincreasing function $C$. (Its value may change from display to display, but it will always be a continuous and nonincreasing function.)

Since $\|G\|_2^2 = \|\bar{S}\|_1/2$, the previous estimate also holds for $\|G\|_2^2$.

4.2.2. *Expliciting the bound.* We deduce from Lemma 2 that

$$2\|S_\tau - \bar{S}\|_1 \leq \|H_\tau + G_\tau - 2G\|_2(4\|G\|_2 + \|H_\tau + G_\tau - 2G\|_2)$$
$$+ \|H_\tau - G_\tau\|_2^2.$$

Making use of Lemma 3 and the bound we obtained for $\|\bar{S}\|_1$, we get, for $N$ large enough and $s$ on $[s_0,\infty)$,

$$N^{2/3}\|S_\tau - \bar{S}\|_1 \leq C(s_0)e^{-s/2}(5C(s_0)e^{-s/2}) + C^2(s_0)e^{-s}.$$

Therefore,

$$N^{2/3}\|S_\tau - \bar{S}\|_1 \leq C(s_0)\exp(-s/2),$$

for yet another continuous nonincreasing function $C$. We remark that if the $C$ appearing in Lemma 3 grows polynomially as $s_0 \to -\infty$, the $C$ appearing in the previous display will also do so.

Finally, combining the previous estimate with Lemma 1, we obtain, when $N$ is large enough and $s$ is in $[s_0,\infty)$,

$$N^{2/3}(P(l_{n,N} \leq s) - F_2(s)) \leq C(s_0)e^{-s}e^{2+2\|\bar{S}\|_1},$$

where $\bar{S}$ is seen as an operator on $L^2([s,\infty))$. Using the fact that $\|\bar{S}\|_1$ is a nonincreasing function of $s$, we finally get

$$M(s_0) = C(s_0)e^{2+2\|\bar{S}\|_1(s_0)}.$$

4.2.3. *Qualitative properties of $M$.* The (continuous and nonincreasing) function $C$ appearing in the course of the proof of Fact 2.2.1 is essentially controlled, especially when $s \to -\infty$ by the maximum of $|\text{Ai}'|$ on $[s,\infty)$. This roughly behaves like $|s|^{1/4}$ when $s$ goes to $-\infty$. A careful look at the intermediate steps in the proof of Fact 2.2.1 shows that $C$ can be chosen to grow polynomially when $s_0 \to -\infty$. On the other hand, this function can be taken to be constant on $[s_1, \infty)$. As explained in Section 4.2.2, such remarks carry over to the choice of $C$ made in the definition of $M$.

Since $\|\bar{S}\|_1(s_0)$ is of order $|s_0|^{3/2}$ at $-\infty$ and appears in an exponential in the definition of $M$, we see that this part will be much more important in the growth of $M$ than the function $C(s_0)$. It also appears that $M$ grows when its argument goes to $-\infty$ and does it fast.

Note that the bound for the Hilbert–Schmidt norm of the difference of the operators (Lemma 3) could be a little improved upon by using better estimates than exponential bounds for Ai and Ai$'$ when their argument is negative and large. Indeed, the analysis presented to obtain Fact 2.2.1 actually provides a polynomially growing bounding function $P$ (of $s_0$) for the difference of the functions we are interested in on intervals of the type $[-s_0, s_1]$. In other words, we could replace $C(s_0)\exp(-s/2)$ by $P(s_0)$ on those intervals $P(s_0)$ growing polynomially to $\infty$ as $s_0 \to -\infty$. The speed of growth of the function $\|\bar{S}\|_1(s)$ when $s$ goes to $-\infty$ nevertheless makes such improvement of little importance.

So our bounding function $M(s_0)\exp(-s)$ is a trade-off between two quantities. $M$ grows—theoretically—very fast when $s_0$ tends to $-\infty$. On the other hand, it seems to not be too large for $s$ nonnegative, or $s$ negative but not too large in absolute value. Finally, the form of the bounding function might shed some light on the empirical findings that the Tracy–Widom approximation is particularly good in the upper tail of the distribution.

4.2.4. *Role of $\gamma$ in the problem.* The analysis we presented to obtain Fact 2.2.1 and through the previous two subsections shows that, given $\gamma$ (the limit of $n/N$ as they go to $\infty$), we can find $M_\gamma$ such that Theorem 2 holds (with $M_\gamma$ taking the place of $M$). We present in Appendix A.6 a study of the dependence (with respect to $\gamma$) of the intermediate functions we obtained. It shows that these functions $M_\gamma$ actually satisfy

$$M_\gamma(s) \leq \tilde{M}(s)(1+\gamma^{-1/2}) \leq 2\tilde{M}(s) \leq M(s),$$

since $\gamma \geq 1$.

TABLE 1
*Empirical quality of new centering and scaling sequence: simulations. The leftmost columns display certain quantiles of the Tracy–Widom distribution. The second column gives the corresponding value of its c.d.f. Other columns give the value of the empirical distribution functions obtained from simulations at these quantiles. The centering and scaling sequences are $\tilde{\mu}_{n,N}$ and $\tilde{\sigma}_{n,N}$*

| Quantiles | TW | $1000 \times 10$ | $4000 \times 10$ | $10000 \times 10$ | $5000 \times 30$ | $4000 \times 100$ | 2*SE |
|---|---|---|---|---|---|---|---|
| −3.73 | 0.01 | 0.010 | 0.014 | 0.013 | 0.010 | 0.009 | 0.002 |
| −3.20 | 0.05 | 0.048 | 0.060 | 0.058 | 0.052 | 0.046 | 0.004 |
| −2.90 | 0.10 | 0.100 | 0.113 | 0.113 | 0.103 | 0.095 | 0.006 |
| −2.27 | 0.30 | 0.300 | 0.313 | 0.312 | 0.304 | 0.295 | 0.009 |
| −1.81 | 0.50 | 0.510 | 0.502 | 0.513 | 0.501 | 0.497 | 0.010 |
| −1.33 | 0.70 | 0.707 | 0.701 | 0.704 | 0.700 | 0.696 | 0.009 |
| −0.60 | 0.90 | 0.899 | 0.904 | 0.901 | 0.900 | 0.893 | 0.006 |
| −0.23 | 0.95 | 0.949 | 0.953 | 0.952 | 0.952 | 0.947 | 0.004 |
| 0.48 | 0.99 | 0.990 | 0.991 | 0.990 | 0.990 | 0.989 | 0.002 |

In some sense this hides the role of $\gamma$, which might appear from the previous display to have little influence on the problem. This parameter nevertheless seems to have an important role in determining $N(s_0, \gamma)$, that is, the integer after which our estimates are correct. The fact that, numerically, large values of $\gamma$ seem to yield an approximation of lower quality in simulations is likely to be a manifestation of the role of $\gamma$ in $N(s_0, \gamma)$.

4.3. *Simulations.* Part of the impetus for this study was practical. We were wondering if it was possible to find new centering and scaling sequences that would improve the quality of the Tracy–Widom approximation in "small" dimensions, a crucial need in Statistics and, more generally, for the relevance in applied fields of such approximations. We made some simulations to see how the sequences we found and chose behaved in practice.

Tables 1–3 were constructed (as in [11]) in the following way: we built $n \times N$ matrices $X$, filled with i.i.d. standard complex Gaussian entries. We made 10,000 simulations for each matrix size, extracting the largest eigenvalue of $X^*X$, recentering and rescaling it through $\tilde{\mu}_{n,N}$ and $\tilde{\sigma}_{n,N}$. The quantiles of the Tracy–Widom law $W_2$ are courtesy of Professor Iain M. Johnstone. The columns corresponding to matrix sizes give the value of the empirical distribution function we found for $l_{n,N}$ at the Tracy–Widom quantiles. If the approximation were "perfect," all the columns would be equal to the second column starting from the left.

The conclusion we can draw from these simulations is that the approximation is remarkably good in the upper tail of the distribution (as also remarked in [11]; this is an excellent property to have for one-sided tests in Statistics) and that the new sequences seem to improve the quality of the approximation over the so-far standard choices, especially in the cases where $n/N$ is quite large. Those are among the most interesting from a statistical point of view, as they are often encountered in "neo-classical" settings ($n$ quite large and $N$ moderately large) and they are situations where the rapidly developing theory of large random covariance matrices provide an alternative to the classical theory (see, e.g., [1], Chapters 11 and 13 and, in particular, Theorem 13.5.2) of multivariate statistics.

TABLE 2
*Empirical quality of new centering and scaling sequence: simulations. The data was generated as in Table 1. We experiment with (increasingly) large $\gamma$*

| Quantiles | TW | $5000 \times 50$ | $20000 \times 50$ | $50000 \times 50$ | $200 \times 5$ | $2000 \times 5$ | $20000 \times 5$ | 2*SE |
|---|---|---|---|---|---|---|---|---|
| $-3.73$ | 0.01 | 0.011 | 0.011 | 0.012 | 0.008 | 0.014 | 0.014 | 0.002 |
| $-3.20$ | 0.05 | 0.052 | 0.050 | 0.056 | 0.042 | 0.058 | 0.066 | 0.004 |
| $-2.90$ | 0.10 | 0.101 | 0.105 | 0.108 | 0.089 | 0.109 | 0.121 | 0.006 |
| $-2.27$ | 0.30 | 0.301 | 0.308 | 0.314 | 0.293 | 0.312 | 0.327 | 0.009 |
| $-1.81$ | 0.50 | 0.500 | 0.498 | 0.509 | 0.502 | 0.509 | 0.523 | 0.010 |
| $-1.33$ | 0.70 | 0.706 | 0.700 | 0.703 | 0.703 | 0.714 | 0.711 | 0.009 |
| $-0.60$ | 0.90 | 0.903 | 0.904 | 0.903 | 0.904 | 0.904 | 0.904 | 0.006 |
| $-0.23$ | 0.95 | 0.951 | 0.951 | 0.953 | 0.956 | 0.953 | 0.952 | 0.004 |
| 0.48 | 0.99 | 0.990 | 0.991 | 0.990 | 0.991 | 0.991 | 0.991 | 0.002 |

TABLE 3
*Empirical quality of new centering and scaling sequence: simulations. The data was generated as in Table 1. We look at relatively small matrices and $\gamma$ quite small (1 and 4)*

| Quantiles | TW | $5 \times 5$ | $10 \times 10$ | $100 \times 100$ | $20 \times 5$ | $40 \times 10$ | $400 \times 100$ | 2*SE |
|---|---|---|---|---|---|---|---|---|
| $-3.73$ | 0.01 | 0 | 0.001 | 0.006 | 0.001 | 0.003 | 0.007 | 0.002 |
| $-3.20$ | 0.05 | 0.002 | 0.012 | 0.041 | 0.017 | 0.026 | 0.041 | 0.004 |
| $-2.90$ | 0.10 | 0.012 | 0.042 | 0.084 | 0.049 | 0.065 | 0.089 | 0.006 |
| $-2.27$ | 0.30 | 0.168 | 0.228 | 0.284 | 0.232 | 0.260 | 0.285 | 0.009 |
| $-1.81$ | 0.50 | 0.412 | 0.454 | 0.491 | 0.452 | 0.472 | 0.490 | 0.010 |
| $-1.33$ | 0.70 | 0.669 | 0.682 | 0.695 | 0.684 | 0.691 | 0.688 | 0.009 |
| $-0.60$ | 0.90 | 0.903 | 0.903 | 0.902 | 0.906 | 0.903 | 0.894 | 0.006 |
| $-0.23$ | 0.95 | 0.952 | 0.952 | 0.951 | 0.953 | 0.952 | 0.948 | 0.004 |
| 0.48 | 0.99 | 0.992 | 0.989 | 0.990 | 0.992 | 0.991 | 0.989 | 0.002 |

## APPENDIX

### A.1. Properties of $r_N$, $\tilde{\mu}_{n,N}$, $\tilde{\sigma}_{n,N}$ and $\alpha_{n,N}$.

A.1.1. *Study of $r_N$.* We are interested in asymptotics for $r_N$. Recall from [11] that

$$r_N^2 = \frac{2\pi \kappa_N^{2\kappa_N} e^{2\kappa_N c_0}}{N!(N+\alpha)!}.$$

Now $\kappa_N^{2\kappa_N} e^{2\kappa_N c_0} = e^{-2\kappa_N} n_+^{n_+} N_+^{N_+} = e^{-n_+ - N_+} n_+^{n_+} N_+^{N_+}$. Stirling's formula gives

$$n! = e^{-n} n^{n+1/2} \sqrt{2\pi} \left(1 + \frac{1}{12n} + \frac{1}{288n^2} + \mathrm{O}\left(\frac{1}{n^3}\right)\right).$$

Rewriting $r_N^2$, we get

$$r_N^2 = \frac{\sqrt{2\pi} e^{-n_+} n_+^{n_+}}{n!} \frac{\sqrt{2\pi} e^{-N_+} N_+^{N_+}}{N!}.$$

Let us therefore focus on $g(n) = \frac{\sqrt{2\pi} e^{-n_+} n_+^{n_+}}{n!}$. Applying Stirling's formula, we get

$$g(n) = e^{-1/2} \left(1 + \frac{1}{2n}\right)^{n+0.5} (1 - x_n + x_n^2 + \mathrm{o}(x_n^2)),$$

with $x_n = 1/(12n) + 1/(288n^2)$. So if we go to second order, we have

$$\left(1 + \frac{1}{2n}\right)^{n+0.5} = e^{1/2}\left(1 + \frac{1}{8n} - \frac{7}{128n^2} + \mathrm{O}\left(\frac{1}{n^3}\right)\right),$$

$$\frac{1}{1+x_n} = 1 - \frac{1}{12n} + \frac{1}{288n^2} + \mathrm{O}\left(\frac{1}{n^3}\right).$$

From this we conclude that

$$g(n) = \left(1 + \frac{1}{24n} - \frac{71}{1152n^2} + \mathrm{O}\left(\frac{1}{n^3}\right)\right),$$

and finally,

$$(\mathrm{A.1}) \quad r_N = 1 + \frac{1}{48}\left(\frac{1}{n} + \frac{1}{N}\right) + \frac{1}{2304nN} - \frac{285}{4608}\left[\frac{1}{n^2} + \frac{1}{N^2}\right] + \mathrm{o}\left(\frac{1}{n^2}, \frac{1}{N^2}\right).$$

A.1.2. *Properties of the centering sequence.* The aim of the discussion that follows is to show that, as we announced in the text, we have

$$\tilde{\mu}_{n,N} - \mu_{n-1,N} = \mathrm{O}(1) \quad \text{and} \quad \tilde{\mu}_{n,N} - \mu_{n,N-1} = \mathrm{O}(1).$$

PROOF. Let us first remark that one can write, with $a = \sigma_{n-1,N}^{1/2}$ and $b = \sigma_{n,N-1}^{1/2}$,

$$\tilde{\mu}_{n,N} = \left(\frac{1}{a} + \frac{1}{b}\right)\left(\frac{1}{a\mu_{n-1,N}} + \frac{1}{b\mu_{n,N-1}}\right)^{-1},$$

so that
$$\tilde{\mu}_{n,N} - \mu_{n-1,N} = \frac{a\mu_{n-1,N}(\mu_{n,N-1} - \mu_{n-1,N})}{b\mu_{n,N-1} + a\mu_{n-1,N}}$$

and
$$\tilde{\mu}_{n,N} - \mu_{n,N-1} = \frac{b\mu_{n,N-1}(\mu_{n,N-1} - \mu_{n-1,N})}{b\mu_{n,N-1} + a\mu_{n-1,N}}.$$

Hence, using the fact that $a\mu_{n-1,N}/(b\mu_{n,N-1}) \to 1$, we just need to show that
$$\mu_{n-1,N} - \mu_{n,N-1} = O(1)$$
to have the announced property. This is an easy task, if a little tedious. To carry it out, we study
$$u_{n+\alpha,N+\beta} = (\sqrt{n+\alpha} + \sqrt{N+\beta})^2.$$

Expanding the square and using the Taylor expansion of $(1+x)^{1/2}$ around 0, we have
$$u_{n+\alpha,N+\beta} = n + N + \alpha + \beta + 2\sqrt{nN} + \sqrt{nN}\left(\frac{\alpha}{n} + \frac{\beta}{N}\right)$$
$$+ \frac{\sqrt{nN}}{2}\left(\frac{\alpha\beta}{nN} - \frac{1}{2}\left(\frac{\alpha^2}{n^2} + \frac{\beta^2}{N^2}\right)\right) + O\left(\frac{N^{1/2}}{n^{5/2}}, \frac{n^{1/2}}{N^{5/2}}\right),$$

if we already account for the fact that $n \asymp N$. Since
$$\mu_{n-1,N} = u_{n-0.5,N+0.5} \quad \text{and} \quad \mu_{n,N-1} = u_{n+0.5,N-0.5},$$

we finally obtain
$$\mu_{n-1,N} - \mu_{n,N-1} = \left(\sqrt{\frac{n}{N}} - \sqrt{\frac{N}{n}}\right) + O\left(\frac{N^{1/2}}{n^{5/2}}, \frac{n^{1/2}}{N^{5/2}}\right).$$

Since $n \asymp N$, $(\sqrt{\frac{n}{N}} - \sqrt{\frac{N}{n}})$ is bounded and we have shown that
$$\mu_{n-1,N} - \mu_{n,N-1} = O(1). \qquad \square$$

A.1.3. *Properties of the scaling sequence.* We are going to show that

(Behavior $\tilde{\sigma}_{n,N}$)
$$\frac{\tilde{\sigma}_{n,N}}{\sigma_{n-1,N}} - 1 = O\left(\frac{1}{n}, \frac{1}{N}\right) \quad \text{and}$$
$$\frac{\tilde{\sigma}_{n,N}}{\sigma_{n,N-1}} - 1 = O\left(\frac{1}{n}, \frac{1}{N}\right).$$

PROOF. We remark that one can formally write that, with obvious substitutions,

$$\tilde{\sigma}_{n,N} = c = (1+\gamma)\left(\frac{1}{a} + \frac{\gamma}{b}\right)^{-1}.$$

Note that in $a, b$ and $\gamma$ are nonnegative in the equation defining $\tilde{\sigma}_{n,N}$ in Theorem 2. Simple algebra leads to $c - a = \gamma a(b-a)/(b+a\gamma)$ and $c - b = b(a-b)/(b+a\gamma)$, from which we see that

$$a \wedge b \leq c \leq a \vee b.$$

Hence, to show the properties stated in the equation (Behavior $\tilde{\sigma}_{n,N}$), it is enough to show that

$$\frac{\sigma_{n-1,N}}{\sigma_{n,N-1}} = 1 + \mathrm{O}\left(\frac{1}{n}, \frac{1}{N}\right).$$

This is of course not a surprise, since when we do a Taylor expansion of all the constituting parts of $\sigma_{n-1,N}$ or $\sigma_{n,N-1}$, the first-order terms are of order $1/n$ or $1/N$. Formally, if we call

$$\tau_{n+\alpha,N+\beta} = \frac{\sqrt{n+\alpha} + \sqrt{N+\beta}}{((n+\alpha)^{-1/2} + (N+\beta)^{-1/2})^{1/3}},$$

we have (as remarked in [11], display preceding equation (A.8))

$$\tau^3_{n+\alpha,N+\beta} = \frac{u^2_{n+\alpha,N+\beta}}{\sqrt{(n+\alpha)(N+\beta)}}.$$

We are interested in ratios of the type $\tau_{n+\alpha,N+\beta}/\tau_{n+\alpha',N+\beta'}$, and so we can focus on

$$\left(\frac{\tau_{n+\alpha,N+\beta}}{\tau_{n+\alpha',N+\beta'}}\right)^3 = \left(1 + \frac{u_{n+\alpha,N+\beta} - u_{n+\alpha',N+\beta'}}{u_{n+\alpha',N+\beta'}}\right)^2 \frac{\sqrt{N+\beta'}\sqrt{n+\alpha'}}{\sqrt{N+\beta}\sqrt{n+\alpha}}.$$

The only case that matters to us is $\alpha = -1/2 = -\beta = -\alpha' = \beta'$. Using the work we did on $\mu_{n-1,N} - \mu_{n,N-1}$ and simple Taylor expansions of $[(1+x)/(1-x)]^{1/2}$, we have easily (after we incorporate $n \asymp N$)

$$\left(\frac{\sigma_{n-1,N}}{\sigma_{n,N-1}}\right)^3 = \left(1 + \frac{\sqrt{n/N} - \sqrt{N/n}}{n+N} + \mathrm{O}\left(\frac{1}{n^2}\right)\right)^2\left(1 + \frac{1}{2}\left(\frac{1}{n} + \frac{1}{N}\right) + \mathrm{O}\left(\frac{1}{n^2}\right)\right).$$

We deduce that indeed

$$\frac{\sigma_{n-1,N}}{\sigma_{n,N-1}} = 1 + \mathrm{O}\left(\frac{1}{n}\right),$$

and the property is shown. □

A.1.4. *Control of first term appearing in $\varphi_\tau$.* The aim of this last part is to show that the $\alpha_{n-1,N}$ introduced in Section 3.4 has the property that

$$\alpha_{n-1,N} = \frac{\sqrt{a_N}\sigma_{n-1,N}^{1/2}\tilde{\sigma}_{n,N}}{\tilde{\mu}_{n,N}} = 1 + \mathrm{O}\left(\frac{1}{N}\right).$$

PROOF. The proof will also show that $\alpha_{n,N-1}$ has the same property. It is very simple:

$$\begin{aligned}
\alpha_{n-1,N} &= \frac{(nN)^{1/4}\sigma_{n-1,N}^{1/2}\tilde{\sigma}_{n,N}}{\tilde{\mu}_{n,N}} \\
&= (nN)^{1/4}\left(\frac{\sigma_{n-1,N}^3}{\mu_{n-1,N}^2}\right)^{1/2}\frac{\tilde{\sigma}_{n,N}}{\sigma_{n-1,N}}\frac{\mu_{n-1,N}}{\tilde{\mu}_{n,N}} \\
&= ((n-1)N)^{1/4}\left(\frac{\sigma_{n-1,N}^3}{\mu_{n-1,N}^2}\right)^{1/2}\frac{\tilde{\sigma}_{n,N}}{\sigma_{n-1,N}}\frac{\mu_{n-1,N}}{\tilde{\mu}_{n,N}}(1-1/n)^{-1/4} \\
&= \frac{\tilde{\sigma}_{n,N}}{\sigma_{n-1,N}}\frac{\mu_{n-1,N}}{\tilde{\mu}_{n,N}}(1-1/n)^{-1/4}.
\end{aligned}$$

Since this is the product of terms which are all $1 + \mathrm{O}(1/N)$, we have

$$\alpha_{n-1,N} = 1 + \mathrm{O}\left(\frac{1}{N}\right). \qquad \square$$

**A.2. A closer look at $\zeta$ and $\hat{f}$.** Recall that we are interested in asymptotics for $\tilde{x}_N(s) = \tilde{\mu}_{n,N} + \tilde{\sigma}_{n,N}s$, where $\tilde{\mu}_{n,N}$ and $\tilde{\sigma}_{n,N}$ are the centering and scaling sequences we chose. So in general, we can write $\tilde{x}_N(s)/\kappa_N = \xi_2 + \tilde{\varepsilon}_N(s)$, and for good choices of $\tilde{\mu}_{n,N}$ and $\tilde{\sigma}_{n,N}$, $\tilde{\varepsilon}_N(s)$ is going to be small provided $s$ is not too big. The rationale for setting up the analysis like this is that, to obtain a finite upper bound for the trace class norm of our operators, we use a "variable split" of $[s_0,\infty)$, involving an interval of the form $[s_1, n^{p/q}s_1]$, so we need a uniform control of $\zeta$ over this type of intervals. [In our applications, $\tilde{\varepsilon}_N(n^{p/q}s_1)$ is small with respect to 1.]

Note that we present asymptotics with $\tilde{\varepsilon}_N(s) \geq 0$. If $\tilde{\varepsilon}_N(s) \leq 0$, the same results hold when we use $2/3(-\zeta)^{3/2} = \int_\xi^{\xi_2}\{-f(t)\}^{1/2}\,dt$. So we can safely apply the estimates that follow without worrying about the sign of $\tilde{\varepsilon}_N(s)$. Finally, let us mention that we will often write $\tilde{\varepsilon}_N$ instead of $\tilde{\varepsilon}_N(s)$ for the sake of readability.

A.2.1. *Asymptotics for $\zeta$.* For the sake of simplicity, we will use the notation $x_N(s)$ (resp. $\varepsilon_N$) instead of $\tilde{x}_N(s)$ (resp. $\tilde{\varepsilon}_N$) in what follows, even though we will be working with $\tilde{\mu}_{n,N}$ and $\tilde{\sigma}_{n,N}$. These centering and scaling sequences are generic sequences in this subsection. They just guarantee that $\varepsilon_N(s) = x_N(s)/\kappa_N - \xi(2)$ is "small."

We have

$$\frac{2}{3}\zeta^{3/2}(x_N(s)) = \int_{\xi_2}^{\xi_2+\varepsilon_N(s)} g_N(t)\,dt$$

$$\text{where } g_N(t) = \frac{(t-\xi_2)^{1/2}(t-\xi_1)^{1/2}}{2t}.$$

Changing variables to $y = (t-\xi_2)/\varepsilon_N$, and denoting $\alpha_N = \xi_2 - \xi_1 = 2\sqrt{4-\omega_N^2}$, $\beta_N = \xi_2 = 2 + \sqrt{4-\omega_N^2}$, we get

$$\frac{2}{3}\zeta^{3/2}(x_N(s)) = \frac{\varepsilon_N^{3/2}\alpha_N^{1/2}}{2\beta_N}\int_0^1 \frac{\sqrt{y}(1+y\varepsilon_N/\alpha_N)^{1/2}}{1+y\varepsilon_N/\beta_N}\,dy.$$

Now since $\alpha_N$ and $\beta_N$ have finite nonzero limits as $N \to \infty$, $y\varepsilon_N/\alpha_N$ and $y\varepsilon_N/\beta_N$ stay small with respect to 1. Hence, we can do a Taylor expansion within the integral and integrate it. We get

$$\frac{2}{3}\zeta^{3/2}(x_N(s)) = \frac{\varepsilon_N^{3/2}\alpha_N^{1/2}}{2\beta_N}\left(\frac{2}{3} + \frac{2}{5}\varepsilon_N\left(\frac{1}{2\alpha_N} - \frac{1}{\beta_N}\right) + O(\varepsilon_N^2)\right).$$

In other words, if we write

$$\eta_N \triangleq 1/(2\alpha_N) - 1/\beta_N,$$

we have

(A.2) $$\zeta^{3/2}(x_N(s)) = \frac{\varepsilon_N^{3/2}\alpha_N^{1/2}}{2\beta_N}\left(1 + \frac{3}{5}\varepsilon_N\eta_N + O(\varepsilon_N^2)\right).$$

A.2.2. *Asymptotics for $\hat{f}$.* First we have $f_N = g_N^2$, so

$$f(\xi_2 + \varepsilon_N) = \frac{\alpha_N\varepsilon_N}{4\beta_N^2}\frac{1+\varepsilon_N/\alpha_N}{(1+\varepsilon_N/\beta_N)^2},$$

from which we deduce that

$$f^{-3/2}(\xi_2 + \varepsilon_N) = \left(\frac{\alpha_N\varepsilon_N}{4\beta_N^2}\right)^{-3/2}[1 - 3\varepsilon_N\eta_N + O(\varepsilon_N^2)].$$

We recall that the quantity we are primarily interested in is $(\kappa_N/\sigma_{n,N}^3)^{1/6}\hat{f}^{-1/4}$, with $\hat{f} = f/\zeta$. Combining the two aforementioned Taylor expansions, and the fact that $\hat{f}^{-1/4} = (\zeta^{3/2}f^{-3/2})^{1/6}$, we get

$$\left(\frac{\kappa_N}{\sigma_{n,N}^3}\right)^{1/6}\hat{f}^{-1/4} = \left(\frac{4\beta_N^2\kappa_N}{\alpha_N\sigma_{n,N}^3}\right)^{1/6}\left(1 - \frac{2}{5}\varepsilon_N\eta_N + O(\varepsilon_N^2)\right).$$

We know from [11], displays (A.7) through (A.8), that $\frac{4\kappa_N^2\beta_N^2}{\kappa_N\alpha_N} = \sigma_{n,N}^3$. Therefore, we conclude that

$$\text{(A.3)} \qquad \left(\frac{\kappa_N}{\sigma_{n,N}^3}\right)^{1/6} \hat{f}^{-1/4} = \left(1 - \frac{2}{5}\varepsilon_N\eta_N + \mathrm{O}(\varepsilon_N^2)\right).$$

A.2.3. *Asymptotics for $\kappa_N^{2/3}\zeta(x_N(s))$*. Using equation (A.2), we get

$$\text{(A.4)} \qquad \begin{aligned}\kappa_N^{2/3}\zeta(x_N(s)) &= \varepsilon_N\kappa_N\left(\frac{\kappa_N\alpha_N}{4\beta_N^2\kappa_N^2}\right)^{1/3}\left(1 + \frac{2}{5}\varepsilon_N\eta_N + \mathrm{O}(\varepsilon_N^2)\right) \\ &= \frac{\varepsilon_N\kappa_N}{\sigma_{n,N}}\left(1 + \frac{2}{5}\varepsilon_N\eta_N + \mathrm{O}(\varepsilon_N^2)\right).\end{aligned}$$

**A.3. A note on the continuity and the variation of $v(\omega,\zeta)$.** We first remind the reader that in order to control the error that the perturbation of the differential equation induces on its solution, Olver ([15], Theorem 11.3.1, page 399) introduced an error-control function: the variation of the function $H$, where $H(x) = \int_0^x y^{-1/2} v(\omega, y)$, where $v$ is the perturbation function of the differential equation. In our case it depends on a parameter, $\omega = 2\lambda/\kappa$, which is not the case in [15]. The issue at stake is as follows: is the error-control function uniformly bounded on the interval where $\omega$ varies? Olver shows that it is finite point by point, but we need to be sure that it is bounded on the interval $[0, 2-\delta]$, $\delta > 0$, to make our rate estimates work: this is what allows us to neglect the term $\varepsilon_2$ that appeared in the analysis, because if the error control function is indeed bounded, then $\varepsilon_2$ is $\mathrm{O}(1/N)$. To explain the strategy, we need to jump a little bit ahead of the proof. Our aim is to find an interval $I_\omega$ (for $\zeta$), possibly depending on $\omega$, where we will be guaranteed to be while doing our asymptotic developments, and for which we have

$$\int_{I_\omega} |u|^{-1/2} |v(\omega, u)|\, du < \chi,$$

where $\chi$ does not depend on $\omega$.

Before proceeding, recall that

$$f(\omega, \xi) = \frac{\xi^2 - 4\xi + \omega^2}{4\xi^2},$$

$$g(\omega, \xi) = -\frac{1}{4\xi^2}.$$

We warn the reader that we drop the dependence of $f$ on $\omega$ temporarily. Also,

$$\text{(E1)} \qquad v(\omega, \zeta) = \frac{1}{\hat{f}^{1/4}}\frac{d^2(\hat{f}^{1/4})}{d\zeta^2} + \frac{g}{\hat{f}} = -\frac{1}{\hat{f}^{3/4}}\frac{d^2(\hat{f}^{-1/4})}{d\xi^2} + \frac{g}{\hat{f}}.$$

As mentioned above, the aim of this subsection is to explain why the aforementioned integral of $|y|^{-1/2}|v|$ remains bounded when $\omega^2 \in [0, 4-\delta]$. Note also that our $\omega$ is in $[0, 2-\rho]$.

From [15], Lemma 11.3.1, we know that the integral is finite for any given $\omega < 2$. Our point here is to give a simple argument that justifies why this is also true uniformly on this interval. As a first step toward the resolution of the problem, we want to show that,

(R1) $\qquad$ for $\omega^2 \in [0, 4-\delta]$, $\qquad \int_0^\infty |y|^{-1/2} |v(\omega, y)|\, dy < \infty.$

Recall also that, for $x \geq \xi_2 = 2 + \sqrt{4-\omega^2}$,

$$\tfrac{2}{3}\zeta^{3/2} = \int_{\xi_2}^\xi f(\omega, t)^{1/2}\, dt,$$

so it is clear that $\zeta(\omega, \xi)$ is continuous in these two variables. Moreover, if $\xi \geq 4$, $\zeta(\omega, \xi)$ is an increasing function of $\omega \geq 0$.

Now using, for instance, (E1) and expanding the formula in terms of $g$, $\zeta$ and $f$ and its derivatives, one can see that, at $+\infty$,

$$v(\omega, \zeta) \sim -\tfrac{1}{4}\zeta^{-2}$$

and that we can actually find a $\zeta_1$, independent of $\omega \in [0, 2-\rho]$ but possibly dependent on $\rho$, such that, for $\zeta \geq \zeta_1$,

$$|v(\omega, \zeta)| \leq \zeta^{-2}.$$

So if we can show that $v(\omega, \zeta)$ is continuous on $[0, 2-\rho] \times [0, \infty)$, the result (R1) follows.

Now recall that Olver in [15], Lemma 11.3.1, shows that, at $\omega$ fixed,

$$\hat{f}(\xi) = \frac{f(\omega, \xi)}{\zeta} = (p(\omega, \xi))^2 \left\{\frac{3}{2} q(\omega, \xi)\right\}^{-2/3}$$

$[p(\omega, \xi) = (\xi - \xi_1)^{1/2}/(2\xi)$, and $q(\omega, \xi) = (\xi - \xi_2)^{-3/2} \int_{\xi_2}^\xi (t - \xi_2)^{1/2} p(\omega, t)\, dt]$, is a well-behaved function of $\xi$, and furthermore, that $v(\omega, \zeta)$ is continuous on intervals of interest. The dependence of $p$ and $q$ on $\omega$ comes from the fact that both $\xi_1$ and $\xi_2$ are functions of $\omega$. Note that in the proof of Lemma 11.3.1 in [15], Olver shows that as $\xi \to \xi_2$, $q(\omega, \xi) \to (2/3)p(\omega, \xi_2(\omega))$. Since $\omega$ is bounded away from 2, $p(\omega, \xi_2)$ is bounded away from 0 on $[0, 2-\delta] \times [2, 4]$. As a matter of fact, $\xi_2(\omega) - \xi_1(\omega)$ is bounded away from 0 (uniformly in $\omega$, when $\omega \in [0, 2-\delta]$). Furthermore, on intervals of the type $[s_0, \infty)$, which are the ones we are focusing on for our original rate of convergence problem, the corresponding $\xi$'s remain bounded away from 2, for large enough $N$, because $\xi_2(\omega)$ is uniformly (in $\omega$) bounded away from 2. This is what makes $\hat{f}(\omega, \xi)$ and $1/\hat{f}(\omega, \xi)$ well-behaved functions on $[0, 2-\rho] \times [2, \infty)$. As we just mentioned, the interval $[2, \infty)$ for $\xi$ is more than enough for our purposes, the leftmost our $\xi$'s can go is close to $\inf_\omega \xi_2(\omega)$ and this quantity is bounded away from 2.

With this information about $\hat{f}$, we can show from equation (E1) and along the lines of the arguments given in Lemma 11.3.1 in [15] that $v(\omega, \xi)$ is a continuous function in its two variables. Since $\zeta(\omega, \xi)$ is invertible, continuous and has a well-behaved inverse, we conclude that $v(\omega, \zeta)$ is continuous on $[0, 2 - \rho] \times [0, \zeta_1]$. Hence, we deduce that

$$\exists \chi \text{ s.t. } \forall \omega \in [0, 2 - \rho] \qquad \int_0^\infty |y|^{-1/2} |v(\omega, y)| \, dy < \chi.$$

The problem is not finished because, for $s$ negative, $\zeta$ will cross 0. But we have already seen that the $\zeta$'s appearing in our analyses will never correspond to $\xi$'s that are less than 2. With this in mind, calling $\zeta_2(\omega)$ the $\zeta$ that corresponds to $\xi = 2$, we pick $I_\omega = [\zeta_2(\omega), \infty)$. Let us now show that this is a satisfactory choice.

The continuity arguments we just mentioned can also be used to show that if $\zeta_2(\omega)$ corresponds to $\xi = 2$, we have

$$\exists \chi \text{ s.t. } \forall \omega \in [0, 2 - \rho] \qquad \int_{\zeta_2(\omega)}^0 |y|^{-1/2} |v(\omega, y)| \, dy < \chi,$$

from which we finally deduce what we needed:

$$\exists \chi \text{ s.t. } \forall \omega \in [0, 2 - \rho] \qquad F(\omega) = \int_{\zeta_2(\omega)}^\infty |y|^{-1/2} |v(\omega, y)| \, dy < \chi.$$

**A.4. Proof of Lemma 2.** Recall that $S_\tau = H_\tau G_\tau + G_\tau H_\tau$, $\bar{S} = 2G^2$. The key inequality in what follows is ([9], Lemma IV.7.2, page 67)

$$\|AB\|_1 \leq \|A\|_2 \|B\|_2.$$

Our problem is to control $\|S_\tau - \bar{S}\|_1$. Note that

$$S_\tau = \frac{(H_\tau + G_\tau)^2 - (H_\tau - G_\tau)^2}{2}.$$

Let $A \triangleq H_\tau + G_\tau$ and $B \triangleq 2G$. We have, of course,

$$A^2 - B^2 = \frac{(A+B)(A-B) + (A-B)(A+B)}{2},$$

and hence,

$$\|A^2 - B^2\|_1 \leq \tfrac{2}{2} \|A - B\|_2 \|A + B\|_2.$$

We can therefore conclude that, since $S_\tau - \bar{S} = (A^2 - B^2 - (H_\tau - G_\tau)^2)/2$,

$$2\|S_\tau - \bar{S}\|_1 \leq \|H_\tau + G_\tau - 2G\|_2 \|H_\tau + G_\tau + 2G\|_2 + \|H_\tau - G_\tau\|_2^2,$$

which is the claim made in Lemma 2.

Since we can show that $H_\tau \to G$ at rate at least $1/3$ and $G_\tau \to G$ at rate at least $1/3$ too, it is enough to choose the centering and scaling so that

$$N^{2/3}\|H_\tau + G_\tau - 2G\|_2 \leq C(s_0)\exp(-s/2),$$

in order to get the $2/3$ rate.

**A.5. Simplified expressions for $\varphi$ and $\psi$ and a remark on $\alpha_N = 0$.** We present in this subsection the derivation of the simplified expressions for $\varphi$ and $\psi$. These are simple algebraic manipulations, but they greatly simplify the rate work. We also explain why the case $\alpha_N = 0$ does not create a specific problem.

A.5.1. *Case of $\varphi$.* Recall that, by definition, we have

$$\varphi(x) = (-1)^N \sqrt{\frac{a_N}{2}} \{\sqrt{N+\alpha}\,\xi_N(x) - \sqrt{N}\,\xi_{N-1}(x)\},$$

with $\xi_k(x) = \varphi_k(x;\alpha)/x$ and $\varphi_k(x;\alpha) = \sqrt{k!/(k+\alpha)!}\,x^{\alpha/2}e^{-x/2}L_k^\alpha(x)$. $L_k^\alpha$ is of course a Laguerre polynomial, as defined, for instance, in [19], Chapter 5. Rewriting the previous expression, we get that, if we call $v_N = (-1)^N\sqrt{a_N/2}$,

$$v_N x^{\alpha/2-1}e^{-x/2}\left\{\sqrt{N+\alpha}\sqrt{\frac{N!}{(N+\alpha)!}}L_N^\alpha(x) - \sqrt{N}\sqrt{\frac{(N-1)!}{(N-1+\alpha)!}}L_{N-1}^\alpha(x)\right\}$$

$$= v_N x^{\alpha/2-1}e^{-x/2}\left\{\sqrt{\frac{N!}{(N+\alpha-1)!}}[L_N^\alpha - L_{N-1}^\alpha]\right\} = \varphi(x).$$

Using formula (5.1.13) in [19], we have $L_N^\alpha - L_{N-1}^\alpha = L_N^{\alpha-1}$. Hence,

$$x^{\alpha/2-1}e^{-x/2}[L_N^\alpha - L_{N-1}^\alpha] = x^{(\alpha-1)/2}e^{-x/2}L_N^{(\alpha-1)}x^{-1/2}.$$

We can therefore conclude that

$$\varphi(x) = (-1)^N\sqrt{\frac{a_N}{2}}\varphi_N(x;\alpha-1)x^{-1/2}.$$

A.5.2. *Case of $\psi$.* In this case, we have, by definition,

$$\psi(x) = (-1)^N\sqrt{\frac{a_N}{2}}\{\sqrt{N}\,\xi_N(x) - \sqrt{N+\alpha}\,\xi_{N-1}(x)\}$$

and after expanding $\xi_N$ and $\xi_{N-1}$ just as before, we get

$$\psi(x) = (-1)^N\sqrt{\frac{a_N}{2}}\frac{x^{\alpha/2}e^{-x/2}}{x}\sqrt{\frac{(N-1)!}{(N+\alpha)!}}\{NL_N^\alpha - (N+\alpha)L_{N-1}^\alpha\}.$$

According to [19], formula (5.1.14), $\{NL_N^\alpha - (N+\alpha)L_{N-1}^\alpha\}/x = -L_{N-1}^{(\alpha+1)}$, from which we deduce that

$$\psi(x) = (-1)^{(N-1)}\sqrt{\frac{a_N}{2}}x^{(\alpha+1)/2}e^{-x/2}\sqrt{\frac{(N-1)!}{(N-1+\alpha+1)!}}L_{N-1}^{\alpha+1}x^{-1/2}$$

$$= (-1)^{(N-1)}\sqrt{\frac{a_N}{2}}\varphi_{N-1}(x;\alpha+1)x^{-1/2}.$$

A.5.3. *Case $\alpha_N = 0$.* This situation might seem a little problematic since $\alpha_N - 1$ appears in the definition of $\varphi$. Note, nevertheless, that going back to the definition of $\varphi$ and $\psi$ (before manipulations), we have $\varphi = \psi$ when $\alpha_N = 0$ (or $n = N$), and both are well defined in terms of Laguerre polynomials $L_N^0$ and $L_{N-1}^0$. So by using the expression obtained for $\psi$ for both $\psi$ and $\varphi$ in that case, we realize that the case $\alpha_N = 0$ does not create any further complications.

**A.6. Dependence of error bounds and $s_1$ on $\gamma$.** First a notational point: we now make explicit the dependence on $\gamma$ of the intermediate bounding functions we obtained by subscripting them with $\gamma$. We will give a sketch of the proof, outlining the reasons for which we can bound $M_\gamma$ by $M(1+\gamma^{-1/2})$, as claimed in Section 4.2. The work is divided into three parts, and everything refers to Section 3.5. We first control $\tilde{\Delta}_{n,N}^F(\tilde{x}_N(s))$, then turn our attention to quantities of the type $B_{n-1,N}(\tilde{x}_N(s))$ and in the last step explain why we control $\|H_\tau - G_\tau\|_2$. All the work is done on $I_{0,N} \cup I_{1,N}$, since $I_{2,N}$ does not pose any difficult problems. Finally, we focus on large values of $\gamma$ since the case $\gamma$ uniformly bounded leads directly to uniform bounds on $M_\gamma$ (nothing prevents us from choosing $M_\gamma$ to be continuous with respect to $\gamma$).

A.6.1. *About $\tilde{\Delta}_{n,N}^F(\tilde{x}_N(s))$.* With our choice of $\tilde{\mu}_{n,N}$ and $\tilde{\sigma}_{n,N}$, we have $c_N = s_N = 0$ so $\alpha_{n-1,N}\tilde{u}_{n-1,N}(s) + \alpha_{n,N-1}\tilde{u}_{n,N-1}(s) = 0$. So to control

$$(\alpha_{n-1,N}u_{n-1,N}(s) + \alpha_{n,N-1}u_{n,N-1}(s))\text{Ai}'(s)$$

in terms of $\gamma$, we essentially just have to control

$$\left|\frac{2\eta_N\kappa_{n-1,N}}{5\sigma_{n-1,N}}\tilde{\varepsilon}_N^2(s) + \text{O}(\tilde{\varepsilon}_N^3(s)\kappa_{n-1,N}/\sigma_{n-1,N})\right|.$$

Elementary manipulations show that $\kappa_{n-1,N} \asymp (1+\gamma)N$, $\sigma_{n-1,N} \asymp \gamma^{1/2}N^{1/3}$, $\eta_N \leq \sqrt{1+\gamma}$, and $\mu_{n-1,N} - \mu_{n,N-1} \asymp \gamma^{1/2} - \gamma^{-1/2}$. By aggregating all this information, we conclude that

$$\left|\frac{2\eta_N\kappa_{n-1,N}}{5\sigma_{n-1,N}}\tilde{\varepsilon}_N^2(s)\right| \leq \chi N^{-2/3}(|s|^2 \vee 1),$$

where $\chi$ is independent of $\gamma$.

One easily deduces from this that

$$N^{2/3}\left|\frac{2\eta_N\kappa_{n-1,N}}{5\sigma_{n-1,N}}\tilde{\varepsilon}_N^2(s) + \mathrm{O}(\tilde{\varepsilon}_N^3(s)\kappa_{n-1,N}/\sigma_{n-1,N})\right||\mathrm{Ai}'(s)|$$
(A.5)
$$\leq C(s_0)e^{-s/2}.$$

On the other hand, the end of Section 3.5 makes very clear that all we have to do in order to bound $R_N(s)$ and $R_{N-1}(s)$ is understand $|u_{n-1,N}(s)|^2$. Since we just shown how to control $u_{n-1,N}(s) - \tilde{u}_{n-1,N}(s)$, we just need to bound $\tilde{u}_{n-1,N}(s)$. By definition, we have

$$|\tilde{u}_{n-1,N}(s)| = \left|\frac{\kappa_{n-1,N}}{\sigma_{n-1,N}}\tilde{\varepsilon}_N(s) - s\right|.$$

Using the aforementioned asymptotic properties of $\sigma_{n-1,N}$ and $\mu_{n-1,N} - \mu_{n,N-1}$, we easily obtain our final estimate:

(A.6) $$|\tilde{u}_{n-1,N}(s)| \leq \chi N^{-1/3}(|s| \vee 1)\frac{\gamma - 1}{\gamma} \leq \chi N^{-1/3}(|s| \vee 1).$$

Combining equations (A.6) and (A.5), we finally get

$$N^{2/3}\tilde{\Delta}_{n,N}^F(\tilde{x}_N(s)) \leq \chi C(s_0)e^{-s/2},$$

a bound that is independent of $\gamma$.

A.6.2. *About $B_{n-1,N}(\tilde{x}_N(s))$.* The key to having a good handle of this quantity is of course to be found in Section 3.4.1. We see that we just need to control $\eta_N\tilde{\varepsilon}_N(s) + s\tilde{\sigma}_{n,N}/\tilde{\mu}_{n,N}$. Using the same estimates as before, we easily get

$$\left|\eta_N\tilde{\varepsilon}_N(s) + s\frac{\tilde{\sigma}_{n,N}}{\tilde{\mu}_{n,N}}\right| \leq N^{-2/3}\chi(|s| \vee 1)\gamma^{-1/2},$$

which implies that

$$N^{2/3}B_{n-1,N}(\tilde{x}_N(s)) \leq C(s_0)e^{-s/2}\gamma^{-1/2}.$$

We can therefore conclude that

$$N^{2/3}|\varphi_\tau(s) + \psi_\tau(s) - \sqrt{2}\mathrm{Ai}(s)| \leq C(s_0)e^{-s/2}(1 + \gamma^{-1/2})$$

and we have now made explicit how $C_\gamma$ depends on $\gamma$ in equation (2).

A.6.3. *About $\varphi_\tau - \psi_\tau$.* Working on $|\varphi_\tau(s) - \text{Ai}(s)/\sqrt{2}|$ is now quite simple since we just have to re-use the estimates we just obtained. We conclude that the dependence of $C_\gamma$ on $\gamma$ is the same in equation (3) as it was in equation (2). We also get the same result for $|\psi_\tau(s) - \text{Ai}(s)/\sqrt{2}|$.

The combination of these three results imply that the bounding function in Theorem 2 has the property

$$M_\gamma(s) \leq (1 + \gamma^{-1/2})M(s).$$

Since $\gamma \geq 1$, our bounding function can be chosen to be independent of $\gamma$.

A.6.4. *About $s_1(\gamma)$.* We mentioned in the course of the proof of Fact 2.2.1 that we could choose $s_1$ independently of $\gamma$. Recall that $s_1$ is chosen as in [11], A.8, page 325. It is defined there as, at fixed $\gamma$, $s_1(\gamma) = c(\gamma)(1 + \delta)$, with $\delta > 0$. Recall that

$$c(\gamma) = \lim_{N \to \infty} \frac{4\xi_2^2 \kappa_N}{\sigma_{n,N}^3(\xi_2 - \xi_1)} \sim_{\gamma \to \infty} \frac{16\gamma/2}{\gamma^{3/2} 4\gamma^{-1/2}}.$$

Since $c(\gamma)$ is a continuous function of $\gamma$ on $[1, \infty)$ having a limit at $\infty$, it is bounded. Hence, the same $s_1$ can be chosen for all $\gamma$'s.

**Acknowledgments.** I would like to thank Professor Iain Johnstone for numerous discussions, preprints, help and, last but not least, for telling me that he thought, based on nonrigorous arguments, that rate 2/3 was achievable. These arguments played a key role in my understanding of the problem, and pointed to the crucial fact that one should "trade-off between operators" to get the higher rate 2/3.

I am also grateful to Professor David Donoho for his comments and support. Many thanks also go to Professors Donald St. P. Richards, Harold Widom and Persi Diaconis for references, correspondence and advice. Finally, I would like to thank an anonymous referee for constructive comments that led to significant improvement of the manuscript.

DEPARTMENT OF STATISTICS
UNIVERSITY OF CALIFORNIA, BERKELEY
367 EVANS HALL
BERKELEY, CALIFORNIA 94720
USA
E-MAIL: nkaroui@stat.berkeley.edu